
\input amstex



\def\spaces{\space\space\space\space\space\space\space\space\space\space}
\def\spacess{\message{\spaces\spaces\spaces\spaces\spaces\spaces\spaces}}
\spacess
\spacess
\message{Annals of Mathematics Style: Current Version: 1.1. June 10, 1992}
\spacess
\spacess

\catcode`\@=11

\hyphenation{acad-e-my acad-e-mies af-ter-thought anom-aly anom-alies
an-ti-deriv-a-tive an-tin-o-my an-tin-o-mies apoth-e-o-ses apoth-e-o-sis
ap-pen-dix ar-che-typ-al as-sign-a-ble as-sist-ant-ship as-ymp-tot-ic
asyn-chro-nous at-trib-uted at-trib-ut-able bank-rupt bank-rupt-cy
bi-dif-fer-en-tial blue-print busier busiest cat-a-stroph-ic
cat-a-stroph-i-cally con-gress cross-hatched data-base de-fin-i-tive
de-riv-a-tive dis-trib-ute dri-ver dri-vers eco-nom-ics econ-o-mist
elit-ist equi-vari-ant ex-quis-ite ex-tra-or-di-nary flow-chart
for-mi-da-ble forth-right friv-o-lous ge-o-des-ic ge-o-det-ic geo-met-ric
griev-ance griev-ous griev-ous-ly hexa-dec-i-mal ho-lo-no-my ho-mo-thetic
ideals idio-syn-crasy in-fin-ite-ly in-fin-i-tes-i-mal ir-rev-o-ca-ble
key-stroke lam-en-ta-ble light-weight mal-a-prop-ism man-u-script
mar-gin-al meta-bol-ic me-tab-o-lism meta-lan-guage me-trop-o-lis
met-ro-pol-i-tan mi-nut-est mol-e-cule mono-chrome mono-pole mo-nop-oly
mono-spline mo-not-o-nous mul-ti-fac-eted mul-ti-plic-able non-euclid-ean
non-iso-mor-phic non-smooth par-a-digm par-a-bol-ic pa-rab-o-loid
pa-ram-e-trize para-mount pen-ta-gon phe-nom-e-non post-script pre-am-ble
pro-ce-dur-al pro-hib-i-tive pro-hib-i-tive-ly pseu-do-dif-fer-en-tial
pseu-do-fi-nite pseu-do-nym qua-drat-ics quad-ra-ture qua-si-smooth
qua-si-sta-tion-ary qua-si-tri-an-gu-lar quin-tes-sence quin-tes-sen-tial
re-arrange-ment rec-tan-gle ret-ri-bu-tion retro-fit retro-fit-ted
right-eous right-eous-ness ro-bot ro-bot-ics sched-ul-ing se-mes-ter
semi-def-i-nite semi-ho-mo-thet-ic set-up se-vere-ly side-step sov-er-eign
spe-cious spher-oid spher-oid-al star-tling star-tling-ly
sta-tis-tics sto-chas-tic straight-est strange-ness strat-a-gem strong-hold
sum-ma-ble symp-to-matic syn-chro-nous topo-graph-i-cal tra-vers-a-ble
tra-ver-sal tra-ver-sals treach-ery turn-around un-at-tached un-err-ing-ly
white-space wide-spread wing-spread wretch-ed wretch-ed-ly Brown-ian
Eng-lish Euler-ian Feb-ru-ary Gauss-ian Grothen-dieck Hamil-ton-ian
Her-mit-ian Jan-u-ary Japan-ese Kor-te-weg Le-gendre Lip-schitz
Lip-schitz-ian Mar-kov-ian Noe-ther-ian No-vem-ber Rie-mann-ian
Schwarz-schild Sep-tem-ber Za-mo-lod-chi-kov Knizh-nik quan-tum Op-dam
Mac-do-nald Ca-lo-ge-ro Su-ther-land Mo-ser Ol-sha-net-sky  Pe-re-lo-mov
in-de-pen-dent ope-ra-tors cyclo-tomoc 
}

\Invalid@\nofrills
\Invalid@\usualspace
\newif\ifnofrills@
\def\nofrills@#1#2{\relaxnext@
  \DN@{\ifx\next\nofrills
    \nofrills@true\let#2\relax\DN@\nofrills{\nextii@}%
  \else
    \nofrills@false\def#2{#1}\let\next@\nextii@\fi
\next@}}
\def\usualspace@#1{\ifnofrills@\def\usualspace{#1}\fi}
\def\addto#1#2{\csname \expandafter\eat@\string#1@\endcsname
  \expandafter{\the\csname \expandafter\eat@\string#1@\endcsname#2}}
\newdimen\bigsize@
\def\big@#1#2{{\hbox{$\left#2\vcenter to#1\bigsize@{}%
  \right.\nulldelimiterspace\z@\m@th$}}}
\def\big{\big@\@ne}
\def\Big{\big@{1.5}}
\def\bigg{\big@\tw@}
\def\Bigg{\big@{2.5}}
\def\raggedcenter@{\leftskip\z@ plus.4\hsize \rightskip\leftskip
 \parfillskip\z@ \parindent\z@ \spaceskip.3333em \xspaceskip.5em
 \pretolerance9999\tolerance9999 \exhyphenpenalty\@M
 \hyphenpenalty\@M \let\\\linebreak}
\def\upperspecialchars{\def\ss{SS}\let\i=I\let\j=J\let\ae\AE\let\oe\OE
  \let\o\O\let\aa\AA\let\l\L}
\def\uppercasetext@#1{%
  {\spaceskip1.2\fontdimen2\the\font plus1.2\fontdimen3\the\font
   \upperspecialchars\uctext@#1$\m@th\aftergroup\eat@$}}
\def\uctext@#1$#2${\endash@#1-\endash@$#2$\uctext@}
\def\endash@#1-#2\endash@{%
\uppercase{#1}\if\notempty{#2}--\endash@#2\endash@\fi}
\def\runaway@#1{\DN@{#1}\ifx\envir@\next@
  \Err@{You seem to have a missing or misspelled \string\end#1 ...}%
  \let\envir@\empty\fi}
\newif\iftemp@
\def\notempty#1{TT\fi\def\test@{#1}\ifx\test@\empty\temp@false
  \else\temp@true\fi \iftemp@}

\font@\tensmc=cmcsc10
\font@\sevenex=cmex7 
\font@\sevenit=cmti7
\font@\eightrm=cmr8 
\font@\sixrm=cmr6 
\font@\eighti=cmmi8     \skewchar\eighti='177 
\font@\sixi=cmmi6       \skewchar\sixi='177   
\font@\eightsy=cmsy8    \skewchar\eightsy='60 
\font@\sixsy=cmsy6      \skewchar\sixsy='60   
\font@\eightex=cmex8 %
\font@\eightbf=cmbx8 
\font@\sixbf=cmbx6   
\font@\eightit=cmti8 
\font@\eightsl=cmsl8 
\font@\eightsmc=cmcsc10 
\font@\eighttt=cmtt8 

\loadmsam
\loadmsbm
\loadeufm
\UseAMSsymbols

\def\penaltyandskip@#1#2{\relax\ifdim\lastskip<#2\relax\removelastskip
      \ifnum#1=\z@\else\penalty@#1\relax\fi\vskip#2%
  \else\ifnum#1=\z@\else\penalty@#1\relax\fi\fi}
\def\nobreak{\penalty\@M
  \ifvmode\def\penalty@{\let\penalty@\penalty\count@@@}%
  \everypar{\let\penalty@\penalty\everypar{}}\fi}
\let\penalty@\penalty

\def\block{\RIfMIfI@\nondmatherr@\block\fi
       \else\ifvmode\vskip\abovedisplayskip\noindent\fi
        $$\def\endblock{\par\egroup$$}\fi
  \vbox\bgroup\advance\hsize-2\indenti\noindent}
\def\endblock{\par\egroup}

\def\logo@{\baselineskip2pc \hbox to\hsize{\hfil\eightpoint Typeset by
 \AmSTeX}}




\font\elevensc=cmcsc10 scaled\magstephalf
\font\tensc=cmcsc10

\font\eightsc=cmcsc10 scaled800

\font\elevenrm=cmr10 scaled \magstephalf
\font\ninerm=cmr9
\font\eightrm=cmr8
\font\sixrm=cmr6
\font\fiverm=cmr5

\font\eleveni=cmmi10 scaled\magstephalf
\font\ninei=cmmi9
\font\eighti=cmmi8
\font\sixi=cmmi6
\font\fivei=cmmi5
\skewchar\ninei='177 \skewchar\eighti='177 \skewchar\sixi='177
\skewchar\eleveni='177

\font\elevensy=cmsy10 scaled\magstephalf
\font\ninesy=cmsy9
\font\eightsy=cmsy8
\font\sixsy=cmsy6
\font\fivesy=cmsy5
\skewchar\ninesy='60 \skewchar\eightsy='60 \skewchar\sixsy='60
\skewchar\elevensy'60

\font\eighteenbf=cmbx10 scaled\magstep3

\font\twelvebf=cmbx10 scaled \magstep1
\font\elevenbf=cmbx10 scaled \magstephalf
\font\tenbf=cmbx10
\font\ninebf=cmbx9
\font\eightbf=cmbx8
\font\sixbf=cmbx6
\font\fivebf=cmbx5

\font\elevenit=cmti10 scaled\magstephalf
\font\nineit=cmti9
\font\eightit=cmti8

\font\eighteenmib=cmmib10 scaled \magstep3
\font\twelvemib=cmmib10 scaled \magstep1
\font\elevenmib=cmmib10 scaled\magstephalf
\font\tenmib=cmmib10
\font\eightmib=cmmib10 scaled 800 
\font\sixmib=cmmib10 scaled 600

\font\eighteensyb=cmbsy10 scaled \magstep3 
\font\twelvesyb=cmbsy10 scaled \magstep1
\font\elevensyb=cmbsy10 scaled \magstephalf
\font\tensyb=cmbsy10 
\font\eightsyb=cmbsy10 scaled 800
\font\sixsyb=cmbsy10 scaled 600
 
\font\elevenex=cmex10 scaled \magstephalf
\font\tenex=cmex10     
\font\eighteenex=cmex10 scaled \magstep3


\def\elevenpoint{\def\rm{\fam0\elevenrm}%
  \textfont0=\elevenrm \scriptfont0=\eightrm \scriptscriptfont0=\sixrm
  \textfont1=\eleveni \scriptfont1=\eighti \scriptscriptfont1=\sixi
  \textfont2=\elevensy \scriptfont2=\eightsy \scriptscriptfont2=\sixsy
  \textfont3=\elevenex \scriptfont3=\tenex \scriptscriptfont3=\tenex
  \def\bf{\fam\bffam\elevenbf}%
  \def\it{\fam\itfam\elevenit}%
  \textfont\bffam=\elevenbf \scriptfont\bffam=\eightbf
   \scriptscriptfont\bffam=\sixbf
\normalbaselineskip=13.95pt
  \setbox\strutbox=\hbox{\vrule height9.5pt depth4.4pt width0pt\relax}%
  \normalbaselines\rm}

\elevenpoint 

\def\ninepoint{\def\rm{\fam0\ninerm}%
  \textfont0=\ninerm \scriptfont0=\sixrm \scriptscriptfont0=\fiverm
  \textfont1=\ninei \scriptfont1=\sixi \scriptscriptfont1=\fivei
  \textfont2=\ninesy \scriptfont2=\sixsy \scriptscriptfont2=\fivesy
  \textfont3=\tenex \scriptfont3=\tenex \scriptscriptfont3=\tenex
  \def\it{\fam\itfam\nineit}%
  \textfont\itfam=\nineit
  \def\bf{\fam\bffam\ninebf}%
  \textfont\bffam=\ninebf \scriptfont\bffam=\sixbf
   \scriptscriptfont\bffam=\fivebf
\normalbaselineskip=11pt
  \setbox\strutbox=\hbox{\vrule height8pt depth3pt width0pt\relax}%
  \normalbaselines\rm}

\def\eightpoint{\def\rm{\fam0\eightrm}%
  \textfont0=\eightrm \scriptfont0=\sixrm \scriptscriptfont0=\fiverm
  \textfont1=\eighti \scriptfont1=\sixi \scriptscriptfont1=\fivei
  \textfont2=\eightsy \scriptfont2=\sixsy \scriptscriptfont2=\fivesy
  \textfont3=\tenex \scriptfont3=\tenex \scriptscriptfont3=\tenex
  \def\it{\fam\itfam\eightit}%
  \textfont\itfam=\eightit
  \def\bf{\fam\bffam\eightbf}%
  \textfont\bffam=\eightbf \scriptfont\bffam=\sixbf
   \scriptscriptfont\bffam=\fivebf
\normalbaselineskip=12pt
  \setbox\strutbox=\hbox{\vrule height8.5pt depth3.5pt width0pt\relax}%
  \normalbaselines\rm}


\def\eighteenbold{\def\rm{\fam0\eighteenbf}%
  \textfont0=\eighteenbf \scriptfont0=\twelvebf \scriptscriptfont0=\tenbf
  \textfont1=\eighteenmib \scriptfont1=\twelvemib\scriptscriptfont1=\tenmib
  \textfont2=\eighteensyb \scriptfont2=\twelvesyb\scriptscriptfont2=\tensyb
  \textfont3=\eighteenex \scriptfont3=\tenex \scriptscriptfont3=\tenex
  \def\bf{\fam\bffam\eighteenbf}%
  \textfont\bffam=\eighteenbf \scriptfont\bffam=\twelvebf
   \scriptscriptfont\bffam=\tenbf
\normalbaselineskip=20pt
  \setbox\strutbox=\hbox{\vrule height13.5pt depth6.5pt width0pt\relax}%
\everymath {\fam0 }
\everydisplay {\fam0 }
  \normalbaselines\rm}

\def\elevenbold{\def\rm{\fam0\elevenbf}%
  \textfont0=\elevenbf \scriptfont0=\eightbf \scriptscriptfont0=\sixbf
  \textfont1=\elevenmib \scriptfont1=\eightmib \scriptscriptfont1=\sixmib
  \textfont2=\elevensyb \scriptfont2=\eightsyb \scriptscriptfont2=\sixsyb
  \textfont3=\elevenex \scriptfont3=\elevenex \scriptscriptfont3=\elevenex
  \def\bf{\fam\bffam\elevenbf}%
  \textfont\bffam=\elevenbf \scriptfont\bffam=\eightbf
   \scriptscriptfont\bffam=\sixbf
\normalbaselineskip=14pt
  \setbox\strutbox=\hbox{\vrule height10pt depth4pt width0pt\relax}%
\everymath {\fam0 }
\everydisplay {\fam0 }
  \normalbaselines\bf}

\hsize=31pc
\vsize=48pc

\parindent=22pt
\parskip=0pt

\widowpenalty=10000
\clubpenalty=10000

\topskip=12pt 

\skip\footins=20pt
\dimen\footins=3in 

\abovedisplayskip=6.95pt plus3.5pt minus 3pt
\belowdisplayskip=\abovedisplayskip


\voffset=7pt\hoffset= .7in

\newif\iftitle

\def\amheadline{\iftitle%
\hbox to\hsize{\hss\currannalsline\hss}\else\line{\ifodd\pageno
\hfill\thetitle\hfill\llap{\elevenrm\folio}\else\rlap{\elevenrm\folio}
\hfill\theauthors\hfill\fi}\fi}

\headline={\amheadline}
\footline={\global\titlefalse}


\def\annalsline#1#2{\vfill\eject
\ifodd\pageno\else 
\line{\hfill}
\vfill\eject\fi
\global\titletrue
\def\currannalsline{\eightrm 
{\eightbf#1} (#2), \thepages}}

\def\titleheadline#1{\def\one{#1}\ifx\one\empty\else
\def\thetitle{{
\let\\ \relax\eightsc\uppercase{#1}}}\fi}

\newif\ifshort

\let\shorttitle\titleheadline

\def\onpages#1#2{\def\thepages{#1--#2}}

\def\thismuchskip[#1]{\vskip#1pt}
\def\ilook{\ifx\next[ \let\go\thismuchskip\else
\let\go\relax\vskip1pt\fi\go}

\def\institution#1{\def\theinstitutions{\vbox{\baselineskip10pt
\def\and{{\eightrm and }}
\def\\{\futurelet\next\ilook}\eightsc #1}}}
\let\institutions\institution

\newwrite\auxfile

\def\startingpage#1{\def\one{#1}\ifx\one\empty\global\pageno=1\else
\global\pageno=#1\fi
\theoremcount=0 \eqcount=0 \sectioncount=0 
\openin1 \jobname.aux \ifeof1 
\onpages{#1}{???}
\else\closein1 \relax\input \jobname.aux
\onpages{#1}{\lastpage}
\fi\immediate\openout\auxfile=\jobname.aux
}

\def\endarticle{\ifRefsUsed\global\RefsUsedfalse%
\else\vskip21pt\theinstitutions%
\nobreak\vskip8pt
\fi%
\write\auxfile{\string\def\string\lastpage{\the\pageno}}}

\outer\def\bye{\endarticle\par \vfill \supereject \end}

\def\document{\let\fontlist@\relax\let\alloclist@\relax
 \elevenpoint}


\newif\ifacks
\long\def\acknowledgements#1{\def\one{#1}\ifx\one\empty\else
\vskip-\baselineskip%
\global\ackstrue\footnote{\ \unskip}{*#1}\fi}

\def\title#1{\titleheadline{#1}
\vbox to80pt{\vfill
\baselineskip=18pt
\parindent=0pt
\overfullrule=0pt
\hyphenpenalty=10000
\everypar={\hskip\parfillskip\relax}
\hbadness=10000
\def\\ {\vskip1sp}
\eighteenbold#1\vskip1sp}}

\newif\ifauthor

\def\author#1{\vskip11pt
\hbox to\hsize{\hss\tenrm By \tensc#1\ifacks\global\acksfalse*\fi\hss}
\ifshort\else\xdef\theauthors{{\eightsc\uppercase{#1}}}\fi%
\vskip21pt\global\authortrue\everypar={\global\authorfalse\everypar={}}}

\def\twoauthors#1#2{\vskip11pt
\hbox to\hsize{\hss%
\tenrm By \tensc#1 {\tenrm and} #2\ifacks\global\acksfalse*\fi\hss}
\ifshort\else\xdef\theauthors{{\eightsc\uppercase{#1 and #2}}}\fi%
\vskip21pt
\global\authortrue\everypar={\global\authorfalse\everypar={}}}


\newcount\theoremcount
\newcount\sectioncount
\newcount\eqcount

\newif\ifspecialnumon

\def\eqnumber=#1 {\global\eqcount=#1 \global\advance\eqcount by-1\relax}
\def\sectionnumber=#1 {\global\sectioncount=#1 
\global\advance\sectioncount by-1\relax}
\def\proclaimnumber=#1 {\global\theoremcount=#1 
\global\advance\theoremcount by-1\relax}

\newif\ifsection
\newif\ifsubsection

\def\intro{\global\authorfalse
\centerline{\bf Introduction}\everypar={}\vskip6pt}

\def\elevenboldmath#1{$#1$\egroup}
\def\mathbold{\hbox\bgroup\elevenbold\elevenboldmath}

\def\section#1{\global\theoremcount=0
\global\eqcount=0
\ifauthor\global\authorfalse\else%
\vskip18pt plus 18pt minus 6pt\fi%
{\parindent=0pt
\everypar={\hskip\parfillskip}
\def\\ {\vskip1sp}\elevenpoint\bf%
\ifspecialnumon\global\specialnumonfalse$\rm\spnum$%
\gdef\sectnum{$\rm\spnum$}%
\else\interlinepenalty=10000%
\global\advance\sectioncount by1\relax\the\sectioncount%
\gdef\sectnum{\the\sectioncount}%
\fi. \hskip6pt#1
\vrule width0pt depth12pt}
\hskip\parfillskip
\global\sectiontrue%
\everypar={\global\sectionfalse\global\interlinepenalty=0\everypar={}}%
\ignorespaces

}


\newif\ifspequation

\let\eqno\leqno 

\newif\ifineqalignno
\let\saveleqalignno\leqalignno                        
\def\leqalignno{\let\eqnu\Eeqnu\saveleqalignno}

\let\eqalignno\leqalignno

\def\sectandeqnum{%
\ifspecialnumon\global\specialnumonfalse
$\rm\spnum$\gdef\eqnum{{$\rm\spnum$}}\else\global\firstlettertrue
\global\advance\eqcount by1 
\ifappend\applett\else\the\sectioncount\fi.%
\the\eqcount
\xdef\eqnum{\ifappend\applett\else\the\sectioncount\fi.\the\eqcount}\fi}

\def\eqnu{\leqno{\hbox{\elevenrm\ifspequation\else(\fi\sectandeqnum
\ifspequation\global\spequationfalse\else)\fi}}}      

\def\Speqnu{\global\setbox\leqnobox=\hbox{\elevenrm
\ifspequation\else%
(\fi\sectandeqnum\ifspequation\global\spequationfalse\else)\fi}}

\def\Eeqnu{\hbox{\elevenrm
\ifspequation\else%
(\fi\sectandeqnum\ifspequation\global\spequationfalse\else)\fi}}

\newif\iffirstletter
\global\firstlettertrue
\def\eqletter#1{\global\specialnumontrue\iffirstletter\global\firstletterfalse
\global\advance\eqcount by1\fi
\gdef\spnum{\the\sectioncount.\the\eqcount#1}\eqnu}

\newbox\leqnobox
\def\outsideeqnu#1{\global\setbox\leqnobox=\hbox{#1}}

\def\eatone#1{}

\def\dosplit#1#2{\vskip-.5\abovedisplayskip
\setbox0=\hbox{$\let\eqno\outsideeqnu%
\let\eqnu\Speqnu\let\leqno\outsideeqnu#2$}%
\setbox1\vbox{\noindent\hskip\wd\leqnobox\ifdim\wd\leqnobox>0pt\hskip1em\fi%
$\displaystyle#1\mathstrut$\hskip0pt plus1fill\relax
\vskip1pt
\line{\hfill$\let\eqnu\eatone\let\leqno\eatone%
\displaystyle#2\mathstrut$\ifmathqed~~\qed\fi}}%
\copy1
\ifvoid\leqnobox
\else\dimen0=\ht1 \advance\dimen0 by\dp1
\vskip-\dimen0
\vbox to\dimen0{\vfill
\hbox{\unhbox\leqnobox}
\vfill}
\fi}

\everydisplay{\lookforbreak}

\long\def\lookforbreak #1$${\def\mathone{#1}
\expandafter\testforbreak\mathone\splitmath @}

\def\testforbreak#1\splitmath #2@{\def\mathtwo{#2}\ifx\mathtwo\empty%
#1$$%
\ifmathqed\vskip-\belowdisplayskip
\setbox0=\vbox{\let\eqno\relax\let\eqnu\relax$\displaystyle#1$}%
\vskip-\ht0\vskip-3.5pt\hbox to\hsize{\hfill\qed}
\vskip\ht0\vskip3.5pt\fi
\else$$\vskip-\belowdisplayskip
\vbox{\dosplit{#1}{\let\eqno\eatone
\let\splitmath\relax#2}}%
\nobreak\vskip.5\belowdisplayskip
\noindent\ignorespaces\fi}


\newif\ifmathqed



\newcount\linenum
\newcount\colnum

\def\spline{\omit&\multispan{\the\colnum}{\hrulefill}\cr}
\def\colcounter{\ifnum\linenum=1\global\advance\colnum by1\fi}

\def\everyline{\noalign{\global\advance\linenum by1\relax}%
\ifnum\linenum=2\spline\fi}

\def\mtable{\bgroup\offinterlineskip
\everycr={\everyline}\global\linenum=0
\halign\bgroup\vrule height 10pt depth 4pt width0pt
\hfill$##$\hfill\hskip6pt\ifnum\linenum>1
\vrule\fi&&\colcounter\hskip12pt\hfill$##$\hfill\hskip12pt\cr}

\def\endmtable{\crcr\egroup\egroup}




\def\xast{*}
\newcount\intable
\newcount\mathcol
\newcount\savemathcol
\newcount\topmathcol
\newdimen\arrayhspace
\newdimen\arrayvspace

\arrayhspace=8pt 
\arrayvspace=12pt 

\newif\ifdollaron

\def\mathalign#1{\def\arg{#1}\ifx\arg\xast%
\let\go\relax\else\let\go\mathalign%
\global\advance\mathcol by1 %
\global\advance\topmathcol by1 %
\expandafter\def\csname  mathcol\the\mathcol\endcsname{#1}%
\fi\go}

\def\arraypickapart#1]#2*{\if#1c \ifmmode\vcenter\else
\global\dollarontrue$\vcenter\fi\else%
\if#1t\vtop\else\if#1b\vbox\fi\fi\fi\bgroup%
\def\one{#2}}

\def\arraystrut{\vrule height .7\arrayvspace depth .3\arrayvspace width 0pt}

\def\array#1#2*{\def\firstarg{#1}%
\if\firstarg[ \def\two{#2} \expandafter\arraypickapart\two*\else%
\ifmmode\vcenter\else\vbox\fi\bgroup \def\one{#1#2}\fi%
\global\everycr={\noalign{\global\mathcol=\savemathcol\relax}}%
\def\\ {\cr}%
\global\advance\intable by1 %
\ifnum\intable=1 \global\mathcol=0 \savemathcol=0 %
\else \global\advance\mathcol by1 \savemathcol=\mathcol\fi%
\expandafter\mathalign\one*%
\mathcol=\savemathcol %
\halign\bgroup&\hskip.5\arrayhspace\arraystrut%
\global\advance\mathcol by1 \relax%
\expandafter\if\csname mathcol\the\mathcol\endcsname r\hfill\else%
\expandafter\if\csname mathcol\the\mathcol\endcsname c\hfill\fi\fi%
$\displaystyle##$%
\expandafter\if\csname mathcol\the\mathcol\endcsname r\else\hfill\fi\relax%
\hskip.5\arrayhspace\cr}

\def\endarray{\crcr\egroup\egroup%
\global\mathcol=\savemathcol %
\global\advance\intable by -1\relax%
\ifnum\intable=0 %
\ifdollaron\global\dollaronfalse $\fi
\loop\ifnum\topmathcol>0 %
\expandafter\def\csname  mathcol\the\topmathcol\endcsname{}%
\global\advance\topmathcol by-1 \repeat%
\global\everycr={}\fi%
}

\def\big#1{{\hbox{$\left#1\vbox to 10pt{}\right.\n@space$}}}
\def\Big#1{{\hbox{$\left#1\vbox to 13pt{}\right.\n@space$}}}
\def\bigg#1{{\hbox{$\left#1\vbox to 16pt{}\right.\n@space$}}}
\def\Bigg#1{{\hbox{$\left#1\vbox to 19pt{}\right.\n@space$}}}


\def\figcaption#1#2#3{\topinsert
\vskip4pt 
\vbox to#3{\vfill}\vskip1sp
\setbox0=\hbox{\eightsc Figure #1.\hskip12pt\eightpoint #2}
\ifdim\wd0>\hsize
\noindent\eightsc Figure #1.\hskip12pt\eightpoint #2
\else
\centerline{\eightsc Figure #1.\hskip12pt\eightpoint #2}
\fi
\vskip16pt
\endinsert}

\def\wfig#1#2#3{\topinsert
\vskip4pt 
\hbox to\hsize{\hss\vbox{\hrule height .25pt width #3
\hbox to #3{\vrule width .25pt height #2\hfill\vrule width .25pt height #2}
\hrule height.25pt}\hss}
\vskip1sp
\centerline{\eightsc Figure #1}
\vskip16pt
\endinsert}

\def\wfigcaption#1#2#3#4{\topinsert
\vskip4pt 
\hbox to\hsize{\hss\vbox{\hrule height .25pt width #4
\hbox to #4{\vrule width .25pt height #3\hfill\vrule width .25pt height #3}
\hrule height.25pt}\hss}
\vskip1sp
\setbox0=\hbox{\eightsc Figure #1.\hskip12pt\eightpoint\rm #2}
\ifdim\wd0>\hsize
\noindent\eightsc Figure #1.\hskip12pt\eightpoint\rm #2\else
\centerline{\eightsc Figure #1.\hskip12pt\eightpoint\rm #2}\fi
\vskip16pt
\endinsert}

\def\tabcaption#1#2{\vskip6pt
\setbox0=\hbox{\eightsc Table #1.\hskip12pt\eightpoint #2}
\ifdim\wd0>\hsize
\noindent\eightsc Table #1.\hskip12pt\eightpoint #2
\else
\centerline{\eightsc Table #1.\hskip12pt\eightpoint #2}
\fi
\vskip6pt}

\def\endinsert{\egroup\if@mid\dimen@\ht\z@\advance\dimen@\dp\z@ 
\advance\dimen@ 12\p@\advance\dimen@\pagetotal\ifdim\dimen@ >\pagegoal 
\@midfalse\p@gefalse\fi\fi\if@mid\smallskip\box\z@\bigbreak\else
\insert\topins{\penalty 100 \splittopskip\z@skip\splitmaxdepth\maxdimen
\floatingpenalty\z@\ifp@ge\dimen@\dp\z@\vbox to\vsize {\unvbox \z@ 
\kern -\dimen@ }\else\box\z@\nobreak\smallskip\fi}\fi\endgroup}

\def\pagecontents{
\ifvoid\topins \else\iftitle\else 
\unvbox \topins \fi\fi \dimen@ =\dp \@cclv \unvbox 
\@cclv 
\ifvoid\topins\else\iftitle\unvbox\topins\fi\fi
\ifvoid \footins \else \vskip \skip \footins \footnoterule 
\unvbox \footins \fi \ifr@ggedbottom \kern -\dimen@ \vfil \fi}


\newif\ifappend

\def\appendix#1#2{\def\applett{#1}\def\two{#2}%
\global\appendtrue
\global\theoremcount=0
\global\eqcount=0
\vskip18pt plus 18pt
\vbox{\parindent=0pt
\everypar={\hskip\parfillskip}
\def\\ {\vskip1sp}\elevenbold Appendix%
\ifx\applett\empty\gdef\applett{A}\ifx\two\empty\else.\fi%
\else\ #1.\fi\hskip6pt#2\vskip12pt}%
\global\sectiontrue%
\everypar={\global\sectionfalse\everypar={}}\nobreak\ignorespaces}

\newif\ifRefsUsed
\long\def\references{\global\RefsUsedtrue\vskip21pt
\theinstitutions
\global\everypar={}\global\bibnum=0
\vskip20pt\goodbreak\bgroup
\vbox{\centerline{\eightsc References}\vskip6pt}%
\ifdim\maxbibwidth>0pt
\leftskip=\maxbibwidth%
\parindent=-\maxbibwidth%
\else
\leftskip=18pt%
\parindent=-18pt%
\fi
\ninepoint
\frenchspacing
\nobreak\ignorespaces\everypar={\amref}%
}

\def\endreferences{\vskip1sp\egroup\global\everypar={}%
\nobreak\vskip8pt\vbox{\thereceived\therevised}
}

\newcount\bibnum

\def\amref#1 {\global\advance\bibnum by1%
\immediate\write\auxfile{\string\expandafter\string\def\string\csname
\space #1croref\string\endcsname{[\the\bibnum]}}%
\leavevmode\hbox to18pt{\hbox to13.2pt{\hss[\the\bibnum]}\hfill}}

\def\bibline{\hbox to30pt{\hrulefill}\/\/}

\def\name#1{{\eightsc#1}}

\newdimen\maxbibwidth
\def\AuthorRefNames [#1] {%
\immediate\write\auxfile{\string\def\string\cite\string##1{[\string##1]}}

\def\amref{\spamref}
\setbox0=\hbox{[#1] }\global\maxbibwidth=\wd0\relax}

\def\spamref[#1] {\leavevmode\hbox to\maxbibwidth{\hss[#1]\hfill}}


\def\footnoterule{\kern-3pt\hrule width1in height.5pt\kern2.5pt}

\def\footnote#1#2{%
\plainfootnote{#1}{{\eightpoint\normalbaselineskip11pt
\normalbaselines#2}}}

\def\vfootnote#1{%
\insert \footins \bgroup \eightpoint\baselineskip11pt
\interlinepenalty \interfootnotelinepenalty
\splittopskip \ht \strutbox \splitmaxdepth \dp \strutbox \floatingpenalty 
\@MM \leftskip \z@skip \rightskip \z@skip \spaceskip \z@skip 
\xspaceskip \z@skip
{#1}$\,$\footstrut \futurelet \next \fo@t}


\newif\iffirstadded
\newif\ifadded

\def\addedlett{}

\def\alltheoremnums{%
\ifspecialnumon\global\specialnumonfalse
\ifadded\global\addedfalse
\iffirstadded\global\firstaddedfalse
\global\advance\theoremcount by1 \fi
\ifappend\applett\else\the\sectioncount\fi.\the\theoremcount\addedlett%
\xdef\theoremnum{\ifappend\applett\else\the\sectioncount\fi.%
\the\theoremcount\addedlett}%
\else$\rm\spnum$\def\theoremnum{{$\rm\spnum$}}\fi%
\else\global\firstaddedtrue
\global\advance\theoremcount by1 
\ifappend\applett\else\the\sectioncount\fi.\the\theoremcount%
\xdef\theoremnum{\ifappend\applett\else\the\sectioncount\fi.%
\the\theoremcount}\fi}

\def\allcorolnums{%
\ifspecialnumon\global\specialnumonfalse
\ifadded\global\addedfalse
\iffirstadded\global\firstaddedfalse
\global\advance\corolcount by1 \fi
\the\corolcount\addedlett%
\else$\rm\spnum$\def\corolnum{$\rm\spnum$}\fi%
\else\global\advance\corolcount by1 
\the\corolcount\fi}


\newcount\corolcount
\def\xcorol{Corollary}
\def\xtheorem{Theorem}
\def\xmaintheorem{Main Theorem}

\newif\ifthtitle

\let\saverparen)
\let\savelparen(
\def\rmparenl{{\rm(}}
\def\rmparenr{{\rm\/)}}
{
\catcode`(=13
\catcode`)=13
\gdef\makeparensRM{\catcode`(=13\catcode`)=13\let(=\rmparenl%
\let)=\rmparenr%
\everymath{\let(\savelparen%
\let)\saverparen}%
\everydisplay{\let(\savelparen%
\let)\saverparen\lookforbreak}}}

\medskipamount=8pt plus.1\baselineskip minus.05\baselineskip

\def\rmtext#1{\hbox{\rm#1}}

\def\proclaim#1{\vskip-\lastskip
\def\one{#1}\ifx\one\xtheorem\global\corolcount=0\fi
\ifsection\global\sectionfalse\vskip-6pt\fi
\medskip
{\elevensc#1}%
\ifx\one\xmaintheorem\global\corolcount=0
\gdef\theoremnum{Main Theorem}\else%
\ifx\one\xcorol\ 
\alltheoremnums 
\else\ \alltheoremnums\fi\fi%
\ifthtitle\ \global\thtitlefalse{\rm(\thethtitle)}\fi.%
\hskip1em\bgroup\let\text\rmtext\makeparensRM\it\ignorespaces}

\def\nonumproclaim#1{\vskip-\lastskip
\def\one{#1}\ifx\one\xtheorem\global\corolcount=0\fi
\ifsection\global\sectionfalse\vskip-6pt\fi
\medskip
{\elevensc#1}.\ifx\one\xmaintheorem\global\corolcount=0
\gdef\theoremnum{Main Theorem}\fi\hskip.5pc%
\bgroup\it\makeparensRM\ignorespaces}

\def\endproclaim{\egroup\medskip}


\def\xproof{Proof}
\def\xremark{Remark}
\def\xcase{Case}
\def\xsubcase{Subcase}
\def\xconjecture{Conjecture}
\def\xstep{Step}
\def\xof{of}

\def\deconstruct#1 #2 #3 #4 #5 @{\def\one{#1}\def\two{#2}\def\three{#3}%
\def\four{#4}%
\ifx\two\empty #1\else%
\ifx\one\xproof%
\ifx\two\xof%
  \ifx\three\xcorol Proof of Corollary \rm#4\else%
     \ifx\three\xtheorem Proof of Theorem \rm#4\else\xone\fi%
  \fi\fi%
\else\xone\fi\fi.}

\def\pickup#1 {\def\this{#1}%
\ifx\this\xproof\global\let\go\demoproof
\global\let\enddemo\endproof\else
\ifx\this\xremark\global\let\go\demoremark\else
\ifx\this\xcase\global\let\go\demostep\else
\ifx\this\xsubcase\global\let\go\demostep\else
\ifx\this\xconjecture\global\let\go\demostep\else
\ifx\this\xstep\global\let\go\demostep\else
\global\let\go\demoproof\fi\fi\fi\fi\fi\fi}

\newif\ifnonum
\def\demo#1{\vskip-\lastskip
\ifsection\global\sectionfalse\vskip-6pt\fi
\def\one{#1 }\def\two{#1*}%
\setbox0=\hbox{\expandafter\pickup\one}\expandafter\go\two}

\def\numbereddemo#1{\vskip-\lastskip
\ifsection\global\sectionfalse\vskip-6pt\fi
\def\two{#1*}%
\expandafter\demoremark\two}

\def\demoproof#1*{\medskip\def\xone{#1}
{\ignorespaces\it\expandafter\deconstruct\xone {} {} {} {} {} @%
\unskip\hskip6pt}\rm\ignorespaces}

\def\demoremark#1*{\medskip
{\it\ignorespaces#1\/} \ifnonum\global\nonumtrue\else
 \alltheoremnums\unskip.\fi\hskip1pc\rm\ignorespaces}

\def\demostep#1 #2*{\vskip4pt
{\it\ignorespaces#1\/} #2.\hskip1pc\rm\ignorespaces}

\def\enddemo{\medskip}

\def\endproof{\ifmathqed\global\mathqedfalse\medskip\else
\parfillskip=0pt~~\hfill\qed\medskip
\fi\global\parfillskip0pt plus 1fil\relax
\gdef\enddemo{\medskip}}

\def\qed{\vbox{\hrule\hbox{\vrule height6pt\hskip6pt\vrule}\hrule}}


\def\proofbox{\parfillskip=0pt~~\hfill\qed\vskip1sp\parfillskip=
0pt plus 1fil\relax}







\def\stripbs#1#2*{\def\one{#2}}

\def\emptyspace{ }
\def\nextthing{}
\def\newline{***}
\def\eatone#1{ }

\def\lookatspace#1{\ifcat\noexpand#1\ \else%
\gdef\nextthing{}\xdef\next{#1}%
\ifx\next\emptyspace%
\let\nextthing\emptyspace\else\ifx\next\newline%
\gdef\nextthing{\eatone}\fi\fi\fi\egroup\nextthing#1}

{\catcode`\^^M=\active%
\gdef\spacer{\bgroup\catcode`\^^M=\active%
\let^^M=\newline\obeyspaces\lookatspace}}

\def\ref#1{\seeifdefined{#1}\expandafter\csname\one\endcsname\spacer}

\def\cite#1{\expandafter\ifx\csname#1croref\endcsname\relax[??]\else
\csname#1croref\endcsname\fi\spacer}


\def\seeifdefined#1{\expandafter\stripbs\string#1croref*%
\crorefdefining{#1}}

\newif\ifcromessage
\global\cromessagetrue

\def\crorefdefining#1{\ifdefined{\one}{}
{\ifcromessage\global\cromessagefalse%
\message{\spaces\spaces\spaces\spaces\spaces\spaces\spaces}%
\message{<Undefined reference.}%
\message{Please TeX file once more to have accurate cross-references.>}%
\message{\spaces\spaces\spaces\spaces\spaces\spaces\spaces}\fi[??]}}

\def\label#1#2*{\gdef\ctest{#2}%
\xdef\currlabel{\string#1croref}
\expandafter\seeifdefined{#1}%
\ifx\empty\ctest%
\xdef\labelnow{\write\auxfile{\noexpand\def\currlabel{\the\pageno}}}%
\else\xdef\labelnow{\write\auxfile{\noexpand\def\currlabel{#2}}}\fi%
\labelnow}

\def\ifdefined#1#2#3{\expandafter\ifx\csname#1\endcsname\relax%
#3\else#2\fi}




\def\articlecontents{
\vskip20pt\centerline{\bf Table of Contents}\everypar={}\vskip6pt
\bgroup \leftskip=3pc \parindent=-2pc 
\def\item##1{\vskip1sp\indent\hbox to2pc{##1.\hfill}}}

\def\endcontents{\vskip1sp\leftskip=0pt\egroup}

\def\journalcontents{\vfill\eject
\def\currannalsline{\hfill}
\global\titletrue
\vglue3.5pc
\centerline{\tensc\hskip12pt TABLE OF CONTENTS}\everypar={}\vskip30pt
\bgroup \leftskip=34pt \rightskip=-12pt \parindent=-22pt 
  \def\\ {\vskip1sp\noindent}
\def\pagenum##1{\unskip\parfillskip=0pt\dotfill##1\vskip1sp
\parfillskip=0pt plus 1fil\relax}
\def\name##1{{\tensc##1}}}


\institution{}
\onpages{0}{0}
\def\lastpage{???}
\def\thetitle{Title ???}
\def\theauthors{Authors ???}
\def\thereceived{}
\def\therevised{}

\gdef\split{\relaxnext@\ifinany@\let\next\insplit@\else
 \ifmmode\ifinner\def\next{\onlydmatherr@\split}\else
 \let\next\outsplit@\fi\else
 \def\next{\onlydmatherr@\split}\fi\fi\let\eqnu\xspliteqnu\next}

\gdef\align{\relaxnext@\ifingather@\let\next\galign@\else
 \ifmmode\ifinner\def\next{\onlydmatherr@\align}\else
 \let\next\align@\fi\else
 \def\next{\onlydmatherr@\align}\fi\fi\let\eqnu\xspliteqnu\next}

\def\spliteqnu{{\tenrm\sectandeqnum}\relax}

\def\xspliteqnu{\tag\spliteqnu}

\catcode`@=12

%
%
%
%
%

\def\for{\  \hbox{ for } \ }
\def\if{ \ \hbox{ if } \ }

\def\where{\  \hbox{ where } \ }
\def\and{\  \hbox{ and } \ }

\def\equal{\buildrel  def \over =}

\def\om{\omega}

\def\ga{\gamma}
\def\ep{\epsilon}

\def\de{\delta}

\def\ka{\kappa}
\def\si{\sigma}


\def\vph{\varphi}

\def\tga{\tilde{\gamma}}

\def\C{\bold{C}}

\def\Q{\bold{Q}}

\def\R{\bold{R}}
\def\N{\bold{N}}
\def\Z{\bold{Z}}

\def\F{\bold{F}}
\def\one{\bold{1}}

\def\0{\bold{0}}

\def\C{\hbox{\bf C}}


\def\f{\Cal{F}}

\def\l{\Cal{L}}

\def\p{\Cal{P}}

\def\h{\Cal{H}}

\def\i{\Cal{I}}
\def\j{\Cal{J}}

\def\lan{\langle}

\def\ran{\rangle}

\def\sgn{\hbox{sgn}}
\font\germ=eufb10 
\def\goth#1{\hbox{\germ #1}}

\def\HH{\goth{H}}
\def\FF{\goth{F}}

\def\AA{\goth{A}}

\def\HH{\hbox{${\Cal H}$\kern-5.2pt${\Cal H}$}}

\def\#{\sharp}


\document

\annalsline{December}{1999}
\startingpage{1}     

\comment
\nopagenumbers
\headline{\ifnum\pageno=1\hfil\else \rightheadline\fi}
\def\rightheadline{\hfil\eightit 
The Macdonald conjecture
\quad\eightrm\folio}

\voffset=2\baselineskip
\endcomment



\title
{ One-dimensional double Hecke algebras \\  and Gaussians}
\shorttitle{Double Hecke algebras}

\acknowledgements{
Partially supported by NSF grant DMS--9877048}

\author{ Ivan Cherednik}

\institutions{
Math. Dept, University of North Carolina at Chapel Hill,   
 N.C. 27599-3250
\\ Internet: chered\@math.unc.edu
}


\intro 
%
%
%
%
%
\vfil

These notes are about applications of the one-dimensional 
double affine Hecke algebra to $q$-Gauss integrals and 
Gaussian sums. The double affine Hecke algebras were designed for a
somewhat different purpose: to clarify the quantum 
Knizhnik-Zamolodchikov equation. Eventually (through the
Macdonald polynomials) they led to a unification of the Harish-Chandra
transform (the zonal case) and the p-adic spherical transform.
It is not just a unification. The new transform is self-dual in 
contrast to its celebrated predecessors. Actually it is very close 
to the Hankel transform (Bessel functions) 
and its recent generalizations from [O3,D,J] (see also [H]).  
The $q$-Gauss integrals are in the focus of the $q$-theory.

The transfer to the roots of unity and Gaussian sums is quite natural as well.
Quantum groups are the  motivation. We generalize and,
at the same time, simplify the Verlinde algebras,
the reduced categories of representations
of quantum groups at roots of unity.
Another interpretation of the Verlinde algebras
is via the Kac-Moody algebras [KL]
(due to Finkelberg for roots of unity). 
A valuable feature is the projective action of $PSL(2,\Z)$
(cf. [K, Theorem 13.8]). It is the foundation of the double
Hecke algebra  technique
and an extension of the framework  of the theory of metaplectic
representations.

The classification of irreducible spherical unitary representations
of one-dimensional double Hecke algebra matches well and generalizes the 
classical formulas for the Gaussian sums.
The complete version
of this paper (arbitrary root systems) is [C5]. See also
[C1].
There are of course other applications of double affine Hecke algebras.
We will not discuss them here (see [C6]). 

This paper  is based on my Harvard
mini-course  (1999) and lectures at University Roma I
(La Sapienza). I am
thankful to D. Kazhdan, C. De Concini,  and the organizers
of the CIME school for the kind invitation 
and publication of the paper.

\vskip 0.2cm
{\bf Gauss integrals and Gaussian sums.}
The starting point of many mathematical and physical theories 
is the formula:
$$
\eqalignno{
&2\int_{0}^\infty e^{-x^2}x^{2k}\hbox{d}x\ =\ \Gamma (k+1/2), \
\Re k>-1/2. 
&(1)
}
$$
Let us give some examples.
\vskip 0.2cm

{(a)} Its  generalization to the Bessel functions, namely,
the invariance of the Gaussian $e^{-x^2}$
with respect to the Hankel transform, 
is the cornerstone of the Plancherel formula.

{(b)} The following ``perturbation'' for the same $\Re k>-1/2$ 
$$
\eqalignno{
&2\int_{0}^{\infty} 
(e^{x^2}+1)^{-1}x^{2k}\hbox{\, d}x = 
(1-2^{1/2-k})\Gamma (k+1/2)\zeta (k+1/2) &(2)
}
$$
is fundamental in the analytic number theory. 

{(c)} The multi-dimensional extension  due to Mehta  with
the integrand   
$\prod_{1\le i<j\le n}(x_i-x_j)^{2k}$ instead of $ x^{2k}$
gave birth to the theory of matrix models and the Macdonald
theory with various applications in mathematics and physics.  

{(d)} Switching to the roots of unity, the Gauss formula 
$$
\eqalignno{
&\sum_{m=0}^{2N-1} e^{{\pi m^2 \over 2N}i}\ =\ 
(1+i)\sqrt{N},\ \ N\in \N&(3)
}
$$
can be considered as a certain counterpart of (1) at $k=0$.

{(e)} Replacing $x^{2k}$ by $\sinh(x)^{2k},$ we come to
the theory of spherical and hypergeometric functions and
to the Harish-Chandra transform. The transform
of  the Gaussian (although very transcendental)
plays an  important role in the harmonic 
analysis on symmetric spaces.

To employ modern mathematics at full potential,
we need to go from Bessel to  hypergeometric functions.
In contrast to the former, the latter  
can be studied, interpreted and generalized
by a variety of methods in the range from 
representation theory and algebraic geometry to integrable models 
and string theory.
However the straightforward passage  $x^{2k}\to \sinh(x)^{2k}$ creates
problems. The spherical transform is not self-dual anymore,
the formula (1) has no $\sinh$-counterpart, and
the Gaussian looses its Fourier-invariance.

\vskip 0.2cm
{\bf Difference setup.}
It was demonstrated recently that these important
features of the classical Fourier transform are restored for the kernel
$$
\eqalignno{
&\delta(x;q,k)\equal\prod_{j=0}^\infty {(1-q^{j+2x})(1-q^{j-2x})\over
(1-q^{j+k+2x})(1-q^{j+k-2x})},\ 0<q<1, \ k\in \C.
&(4)
}
$$ 
Actually the self-duality of the corresponding transform
can be expected a priori because the Macdonald truncated theta-function 
$\delta$ is 
a unification of
$\sinh(x)^{2k}$ and the Harish-Chandra function ($A_1$) serving the
inverse spherical transform.

As to (1), setting $q=\exp(-1/a),\ a>0,$
$$
\eqalignno{
&(-i)\int_{-\infty i}^{\infty i}q^{-x^2}\delta(x;q,k)\hbox{\, d}x=
{2\sqrt{a\pi}}\prod_{j=0}^\infty
{1-q^{j+k}\over
 1-q^{j+2k}},\ \ \Re k>0.
&(5)
}
$$
Here both sides are well-defined for all $k$ except for the poles 
but coincide only when $ \Re k>0$, worse than in (1). This can be fixed
as follows:  
$$
\eqalignno{
&(-i)\int_{1/4-\infty i}^{1/4+\infty i}q^{-x^2}\mu\hbox{\, d}x=
{\sqrt{a\pi}}\prod_{j=1}^\infty
{1-q^{j+k}\over
 1-q^{j+2k}},\ \Re k>-1/2 \hbox{\ \ for}&(6)\cr
&\mu(x;q,k)\equal \prod_{j=0}^\infty {(1-q^{j+2x})(1-q^{j+1-2x})\over
(1-q^{j+k+2x})(1-q^{j+k+1-2x})},\ 0<q<1, \ k\in \C.
&(7)
}
$$
The limit of (6)  multiplied by $(a/4)^{k-1/2}$
as $a\to \infty$ is (1) in the imaginary variant.

One can make (5) entirely algebraic replacing 
$\ga^{-1}=q^{-x^2}$ by its expansion
$$\tilde{\ga}^{-1}\ =\ \sum_{-\infty}^{+\infty}q^{n^2/4}q^{nx}
$$
and using Const Term$(\sum c_nq ^{nx})=c_0:$
$$
\eqalignno{
&\hbox{Const\ Term\ }(\tilde{\ga}^{-1}\delta)=
2\prod_{j=0}^\infty
{1-q^{j+k}\over
 1-q^{j+2k}}.
&(8)
}
$$

\vskip 0.2cm
{\bf Jackson sums.}
A most promising feature of special  $q$-functions is 
the possibility to replace the integrals by sums,
the Jackson integrals. 

Let $\int_\sharp$ be the integration for the 
 path which  begins at $z=\epsilon i+\infty$, moves
to the left till $\epsilon i$, then down
through  the origin to  $-\epsilon i$, and then
returns down the positive real axis to $-\epsilon i+\infty$
(for small $\epsilon$). Then for $|\Im k|<2\epsilon, \Re k>0,$
$$
\eqalignno{
&{1\over 2i}\int_{\sharp} q^{x^2} \delta\hbox{\, d}x\ =
-{a\pi\over 2}\prod_{j=0}^\infty
{(1-q^{j+k})(1-q^{j-k})\over
 (1-q^{j+2k})(1-q^{j+1})}\, \lan \ga\ran,\cr
&\lan \ga\ran =\sum_{j=0}^\infty q^{{(k-j)^2\over 4}}
{1-q^{j+k}\over
1-q^{k}} \prod_{l=1}^j
{1-q^{l+2k-1}\over
 1-q^{l}}\ =
&(9)\cr
& q^{{k^2\over 4}}\prod_{j=1}^\infty
{(1-q^{j/2})(1-q^{j+k})(1+q^{j/2-1/4+k/2})(1+q^{j/2-1/4-k/2})\over
(1-q^j)}.
}
$$
The sum $\lan \ga\ran$ 
is the Jackson integral of $\ga=q^{x^2}$ for a special choice 
($-k/2$) of the starting point
and $\mu/\mu(-k/2)$ taken as the measure. 
The convergence of the sum (9) is for all $k.$ 
 
When $q=\exp(2\pi i/N)$ and $k$ is a positive integer 
$\le N/2$ we come to the  Gauss-Selberg-type
sums:
$$
\eqalignno{
&\sum_{j=0}^{N-2k} q^{{(k-j)^2\over 4}}
{1-q^{j+k}\over
1-q^{k}} \prod_{l=1}^j
{1-q^{l+2k-1}\over
 1-q^{l}}=
\prod_{j=1}^k
(1-q^{j})^{-1}\sum_{m=0}^{2N-1} q^{m^2/4}.&(10)
}
$$
They resemble the modular Gauss-Selberg sums but the difference 
is dramatical. The latter are calculated in the finite fields and are 
embedded into roots of unity right before the final summation.
Our sums are defined entirely in cyclotomic fields.

Substituting $k=[N/2]$ we  arrive at a generalization of  (3).
Let us consider  $N=2k$ only (the case $N=2k+1$ is very 
similar):
$$
\eqalignno{
&\sum_{m=0}^{2N-1} q^{m^2/4}\ =\ q^{k^2\over 4}\Pi \for
\Pi=\prod_{j=1}^k (1-q^{j}).
}
$$
First, $\Pi\bar{\Pi}=(X^{N}-1)(X+1)(X-1)^{-1}(1)=2N.$
Second, $\arg(1-e^{i\phi})=\phi/2-\pi/2$ when 
$0< \phi< 2\pi,$ and therefore 
$$\arg\Pi = {\pi\over N} {k(k+1)\over 2}-{\pi k\over 2} =
{\pi (1-k)\over 4}.$$

To consider general primitive roots
$q=\exp(2\pi i l/N)$ as $(l,N)=1,$ one needs to control
the set of $\arg(q^{j})$ for $1\le j\le k.$ This leads
to a variant of the quadratic reciprocity. We will skip it.

Finally, 
$$\arg( q^{k^2\over 4}\Pi)={\pi (1-k)\over 4}+
{\pi k\over 4}={\pi\over 4},\and 
 q^{k^2\over 4}\Pi=\sqrt{N}(1+i).$$

%
%
%
\section{Q-Mellin transform}

Conceptually, (5) and (9) are close.
They differ only by the choice of integration and the
sign of the Gaussian. 
However the second formula looks
more involved. Let us clarify this. Following [C3],
we are going to deduce both
formulas from a certain generalization of the key property of
the classical Mellin transform. 

First we will introduce the shift operator.
In the setup of this paper, it
is the so-called Askey-Wilson operator,
the following  $q$-deformation of the differentiation:
$$
\eqalignno{
&Sf(x)\equal (q^x-q^{-x})^{-1}\, (p^{-1}-p)\, f(x),\cr 
&pf(x)=f(x+{1\over2}),\ q=\exp(-1/a)>0.
&\eqnu
\label\shift\eqnum*
}
$$
For instance,
$$
\eqalignno{
&S\ga \ =\ -q^{1/4}\ga,\ \ S\ga^{-1}\ =\ q^{-1/4}\ga^{-1} \for
\ga=q^{x^2}=\exp(-x^2/a).
&\eqnu
\label\shiga\eqnum*
}
$$
The name is ``shift operator''  because its action on the basic
($q$-difference) hypergeometric
function results in the
shift of the parameters [AW]. 

For instance,
$$
\eqalignno{
&Sp^{(k)}_n(x) \ =\ (q^{-n/2}-q^{n/2})\, p^{(k+1)}_{n-1}(x),
\ n=1,2,\cdots
&\eqnu
\label\ship\eqnum*
}
$$
for the {\it Rogers polynomials}  
$\{p_n(x)\} (n=0,1,2,\ldots)$ [AI].  They are $x$-even
polynomials in terms of $q^{mx}$ 
$(m\in \Z),$ pairwise orthogonal with respect to
the pairing $\lan f,g\ran =$ Const Term $(fg\de)$  
(here $\de(x;q,k)$  from (4) is replaced by the
corresponding Laurent series), and normalized by the condition
$$
\eqalignno{
& p^{(k)}_n=q^{nx}+q^{-nx}+\hbox{\ lower\  powers,\ except\ for\ }  
p^{(k)}_0=1.
&\eqnu
\label\progers\eqnum*
}
$$
The basic hypergeometric function makes $n$ arbitrary complex  
and contains one more parameter covering the
$BC$-case.

In (\ref\progers), $k\neq -1,-2,-3,\ldots,-m+1$ modulo $2\pi a i\Z.$
For example,
$$
\eqalignno{
&p^{(k)}_1=q^x+q^{-x},\ p^{(k)}_2=q^{2x}+q^{-2x}+{(1-q^k)(1+q)\over 1-q^{k+1}}.
&\eqnu
\label\ptwo\eqnum*
}
$$

The multi-dimensional generalization of $\{ p\}$ is due
to Macdonald (see e.g. [M1]).
The differential shift operators for arbitrary
root systems were introduced by Opdam (see [O1] and also [He] for
the interpretation via the Dunkl operators). The 
difference ones were considered in my papers. They 
depend on $k$ in contrast to the simplest case considered
here ($A_1$). The celebrated constant term conjecture
was verified using these operators.
They are also used in [O2,C3] for analytic
continuations. 
 
\vskip 0.2cm
{\bf Shift-formula.}
Let $\pi a>\ep>0$ and  $\int$ be one of the following integrations:
$$
\eqalignno{
&\int_{im} =  {1\over 2i}\int_{\ep-\infty i}^{\ep+\infty i}dx,\ \
\int_{\sharp} = {1\over 2i}\int_{\infty+\ep i}^{\infty-\ep i }dx.
&\eqnu
\label\ints\eqnum*
}
$$
The path of the second integration 
begins at $x=-\ep i+\infty$, moves
to the left down the positive real axis till $-\ep i$, then circles the origin
and returns up the positive real axis to $\ep i+\infty$
(for small $\ep>0$).

Given $u>0$,  an even function  $E(x)$ is called $u$-regular,
if it is analytic in
$$
\eqalignno{
&D_{im} = \{x\ \mid\ -u-\ep\le \Re x\le u+\ep\} \for \int_{im},
&\eqnu\cr
\label\doms\eqnum*
&D_{\sharp} = 
\{x\ \mid \ -\ep\le \Im x\le \ep,\ -u\le \Re x\le u\} \for \int_{\sharp},
&\eqnu
\label\domsh\eqnum*
}
$$
and continuous on the integration paths and the  boundary of $D.$
It is also assumed to have
coinciding continuous (in terms of $\ka$)  limits
$$E^{\pm}(\ka+\nu i)\ =\ 
\lim_{\nu\to \pm \infty} E(x)\,\de(x;q,\ka+\nu i)(x)
$$
in the case of the imaginary integration $\int_{im}.$
We introduce the $q$-{\it Mellin transform} of $E(x)$ by the formula
$$
\eqalignno{
&\psi_k(E)= \pi_k^{-1}\int E\de(x;q,k),\ 
\pi_k\equal \prod_{j=0}^\infty{1-q^{k+j}\over 1-q^{2k+j}}.
&\eqnu
\label\mellin\eqnum*
}
$$ 
The following theorem can be directly deduced from (\ref\ship).

\proclaim{Theorem} Given  $\ep>0,$ let $K_\ep$ be the set of $k$ such that
$\Re k>2\ep$ for the imaginary
integration,  and $\{\Re k>0,\ |\Im k|<2\ep\}$ 
for the sharp one
$\int_\sharp.$ Assuming that
$E(x)$ is $1$-regular, $k\in K_\ep,$
and the integrals below are  well-defined,
$$
\eqalignno{
&\psi_k(E)\ =\ (1-q^{k+1})\psi_{k+1}(E)+q^{k+3/2}\psi_{k+2}^{(2)}(E)\cr
&\for \psi_{k}^{(2)}(E)\ =\ \pi_k^{-1}\int S^2(E)\de(x;q,k).
&\eqnu
 \label\shiftwo\eqnum*
}
$$
\label\SHIFTWO\theoremnum*
\endproclaim

\vskip 0.2cm
{\bf Gaussians.}
Let us take $E=\ga^{-1}=q^{-x^2}$ for the imaginary integration provided
that $0<q<1.$ Since $S\ga^{-1}=q^{-1/4}\ga^{-1}$ and $S\ga=-q^{1/4}\ga$ we
come to the following equation for $\psi_k= \psi_k(\ga^{-1})$ 
$$
\eqalignno{
&\psi_k\ =\ (1-q^{k+1})\psi_{k+1}+q^{k+1}\psi_{k+2},\ \Im k>0,
&\eqnu
 \label\psiim\eqnum*
}
$$
which readily results in 
$$
\eqalignno{
&C_{k+1}=C_k \for 
C_k=(-1)^k q^{(k+1)(k+2)\over 2}(\psi_k-\psi_{k+1}).
&\eqnu
 \label\psiimc\eqnum*
}
$$
The function $C_k$ can be extended analytically to all $k\in\C.$
Taking into consideration the zeros of $\C_k$
(at $\Z/2$) and its periodicity properties (pure in the real 
direction and with a multiplicator in the imaginary direction) we conclude
that $C=0.$ This, in its turn, is sufficient to establish (5).

When the integration is sharp and $E$ is analytic,  
the analytic continuation of  $\psi_k(E)$ to all $k$ is given by 
Cauchy's theorem:
$$
\eqalignno{
&\psi_k(E)\ =\ -{a\pi\over 2} \prod_{j=0}^\infty
{1-q^{j-k}\over 1-q^{j+1}}\ \times\cr
&\sum_{j=0}^\infty q^{-kj}
{1-q^{j+k}\over
1-q^{k}} \prod_{l=1}^j
{1-q^{l+2k-1}\over  1-q^{l}} E({k+j\over 2}).
&\eqnu
 \label\psishc\eqnum*
}
$$

Let $E=\ga$ and  $\psi_k=\psi_k(\ga).$
Then
$$
\eqalignno{
&\psi_k = (1-q^{k+1})\psi_{k+1}+q^{k+2}\psi_{k+2} \and C_{k+1}=C_k
\for \cr
&\phi_k+\phi_{k+1}=C_k  q^{(k+1)(k+2)\over 2}\where
\phi_k\equal q^{k(k+1)\over 2}\psi_k.
&\eqnu
 \label\psish\eqnum*
}
$$
Following the imaginary case, 
we conclude that $C=0$ and 
$\phi_k+\phi_{k+1}=0,$ but the end of the proof of (9) is a bit more involved.
We use (\ref\psishc) to check that 
$$\phi_{k+2p}=q^{3kp+p^2}\phi_k \for p=2\pi i a,\ q=\exp(-1/a).$$ 
Also $\phi_k$ has zeros at $\Z_+.$
Eventually we come to (9). 

We note that there is a better proof,
simpler and  more algebraic, which is  parallel to 
the one from [C1] (arbitrary root systems). Also both identities
(5) and (9) can be deduced from  Bailey's ${}_6\Psi_6$ summation
theorem, as was observed by Macdonald and Andrews.  

Actually these formulas are important because they 
are ingredients of the following, hopefully new,
general identities:
$$
\eqalignno{
&\lan p^{(k)}_m(x)\,p^{(k)}_n(x)\,q^{-\varsigma x^2} \ran\ =\cr
&q^{\, \varsigma\,{m^2+n^2+2k(m+n)\over 4}}\, p^{(k)}_m({n+k\over 2})\,
p^{(k)}_n({k\over 2})\,\lan q^{-\varsigma x^2} \ran.
&\eqnu
\label\ggmnid\eqnum*
}
$$
Here $\varsigma=1,-1$ correspond respectively to
$$
\eqalignno{
&\lan f(x)\ran_\circ\equal \hbox{Const\ Term}(f(x)\mu(x;q,k)) \and 
&\eqnu
\label\consto\eqnum*
\cr
&\lan f(x)\ran_\bullet\equal 
&\eqnu
\label\sharone\eqnum*
\cr
&\sum_{j=0}^\infty q^{-kj} {1-q^{j+k}\over 1-q^{k}} \prod_{l=1}^j
{1-q^{l+2k-1}\over  1-q^{l}} f({k+j\over 2}).
}
$$
So we switched from  $\int_{im}\de(x;q,k) dx $ and
$\int_{\sharp}\de(x;q,k) dx$ to more algebraic $\lan\ \ran.$ 
For even $f(x),$ it is just the matter of normalization.
Here $q^{\mp x^2}$ serves respectively $\lan\ \ran_\circ$ and 
$\lan\ \ran_\bullet.$
In the case of $\lan\ \ran_\bullet$
we still assume that $0<q<1,$ which is not necessary anymore for 
$\lan\ \ran_\circ$ after $q^{-x^2}$ is replaced by the Laurent
series.

Recalculating (5) and (9):
$$
\eqalignno{
&\lan q^{-x^2}\ran_\circ\ =\
\prod_{j=1}^\infty
{1-q^{j+k}\over
 1-q^{j+2k}},
&\eqnu
\label\ctga\eqnum*\cr
&\lan q^{x^2}\ran_\bullet\ =\
\prod_{j=1}^\infty
{1-q^{j+k}\over 1-q^j}\sum_{j=-\infty}^{\infty}q^{(k+j)^2/ 4}=
&\eqnu
\label\sharpga\eqnum*
\cr
&q^{{k^2\over 4}}\prod_{j=1}^\infty
{(1-q^{j/2})(1-q^{j+k})(1+q^{j/2-1/4+k/2})(1+q^{j/2-1/4-k/2})\over
(1-q^j)}.
}
$$

The so-called group case $k=1$ is a good exercise.
The Rogers polynomials
become the $SL_2$-characters: 
$$p_n(x)\ =\ {q^{(n+1)x}-q^{-(n+1)x}\over
q^{x}-q^{-x}}.
$$ 
The calculation is almost equally simple for arbitrary root
systems (use Weyl's character formula).

%
%
%
\section{Double Hecke algebras}

Double Hecke algebras provide  justifications and
generalizations. In the $A_1$-case, 
$\HH\equal \C[\Cal{B}_q]/((T-t^{1/2})(T+t^{-1/2}))$ for
the group algebra of the group $\Cal{B}_q$ generated by
$T,X,Y,q^{1/4}$ with the relations
$$
\eqalignno{
&TXT=X^{-1},\ T^{-1}YT^{-1}=Y^{-1},\ Y^{-1}X^{-1}YXT^2=q^{-1/2}
&\eqnu
\label\douba\eqnum*
}
$$
for central $q^{1/4},t^{1/2}$.
 
The starting point of the theory 
is the PBW theorem.  Any element of $H\in\HH$ can be uniquely expressed 
in the form
$$
\eqalignno{
H=\sum 
c_{i\epsilon j}X^i T^\epsilon Y^j\quad (c_{i\epsilon j}\in\C)
\where i,j\in \Z,\epsilon=0,1.
&\eqnu
\label\pbwt\eqnum*
}
$$
Permuting $X,T,Y$ one gets 5 more statements of this kind.

Renormalizing $T\to q^{-1/4}T,\ X\to q^{1/4}X,\
 Y\to q^{-1/4}Y,$ we can assume that $q=1$ in  $\Cal{B}_q.$
However this will change the quadratic relation. The group
$\Cal{B}_1$ is close to the fundamental group of the 
$\{ E\times E\setminus \hbox{diag}\}/\bold{S}_2$ for the elliptic
curve $E,$ which is a special case of the calculation due to  Birman.
Without going into detail, let us mention that 
$T$ is the half-turn about the diagonal,
$X,Y$ correspond to  the  ``periods'' of $E.$ 

The topological interpretation is helpful
to see that the central
extension  $PSL^c_2(\Z)$ of
$PSL_2(\Z)$ (Steinberg)
acts on $\Cal{B}_1$ and $\HH$.
The automorphisms corresponding to
the generators ${11\choose 01},\
{10\choose 11}$
are as follows: 
$$
\eqalignno{
&\tau_+: Y\to q^{-1/4}XY,\ X\to X,\ \ \tau_-:  X\to q^{1/4}YX,\ Y\to Y.
&\eqnu
\label\taupm\eqnum*
}
$$
They fix $T,q,t.$
This of course can be readily checked without topology. 

Formally,  $\tau_+$ is the conjugation by
$q^{x^2}$ for  $X$ represented (here and further) in the form
$X = q^x.$  
This will be used to calculate  the $q$-Fourier transform
of the Gaussian multiplied by the Macdonald polynomials.       

The Fourier transform on the generalized functions can be 
associated 
with the anti-involution 
$$
\eqalignno{
&\vph:\ X\mapsto Y^{-1}\mapsto X,\ T\mapsto T,\ q,t\mapsto q,t
&\eqnu
\label\phiant\eqnum*
}
$$
The existence of $\vph$ can be easily deduced from (\ref\douba) too.

The group
$PSL^c_2(\Z)$ and $\vph$ act on  $\HH$ for all root systems. A direct
algebraic proof  of this  important fact is known (the author and Macdonald).
However the considerations are more involved, especially if the root system is
not of the $A$-type.  
When $t=1$ we get the well-known action of  
$SL_2(\Z)$ on the Weyl and  Heisenberg algebras (the latter as
$q\to 1$). 

\vskip 0.2cm
{\bf Macdonald polynomials}. Generalizing the Rogers polynomials,
we introduce the  Mac\-do\-nald (non\-sym\-met\-ric) 
polynomials as eigenfunctions of  $\widehat{Y}$ in the following 
$\HH$-representation
in the space $\l$ of the  Laurent polynomials of $X=q^x:$
$$
\eqalignno{
&T\mapsto \widehat{T}=
t^{1/2}s+(q^{2x}-1)^{-1}(t^{1/2}-t^{-1/2})(s-1), \ 
Y\mapsto \widehat{Y}= sp\widehat{T}
&\eqnu
\label\tyopers\eqnum*
}
$$
for the reflection $sf(x)=f(-x)$ and the translation
$pf(x)=f(x+1/2).$ It is nothing else but the representation
of $\HH$ induced from the character $\chi(T)=t^{1/2}=\chi(Y)$
of $\h_Y=\lan T,Y\ran .$

Explicitly:
$$
\eqalignno{
&\widehat{Y}e_n^{(k)}\ =\ 
q^{-n_{\sharp}}e_n,\ \ n_{\sharp}\equal {n+\sgn(n-1/2)k\over 2},\cr
&e_n^{(k)}=\sum_{m=-n}^{n} c_n^m q^{mx},\ c_n^n=1\hbox {\ (normalization),\ }
\ c_n^{-n}=0 \for n>0.
&\eqnu
\label\epspol\eqnum*
}
$$
From now on we set $t=q^k.$
Here we need to take $k\not\in -\N/2$ to ensure that the spectrum of
$Y$ in $\l$ is simple (i.e. $n_\sharp\neq m_\sharp$ for $n\neq m$). 
In this section  $q$ is generic.

The connection with the Rogers polynomials is via the symmetrization:
$$
\eqalignno{
&p_n\ =\ (1+t^{1/2}T)\,e_{n}\ = \ (1+s)\,
\bigl({t-q^{2x}\over 1-q^{2x}}\, e_n \bigr) \for n>0.
&\eqnu
\label\epsro\eqnum*
}
$$

Let us renormalize $\{e\}.$ We call $\ep_n^{(k)}=e_n^{(k)}/e_n^{(k)}(-k/2)$
{\it spherical polynomials}. The motivation 
is the following duality $\ep_n(m_\sharp)=\ep_m(n_\sharp)$ 
for all $n,m.$ Upon the symmetrization, 
$p_n(m_\sharp)p_m(-k/2)=p_m(n_\sharp)p_n(-k/2)$ for 
$n,m\ge 0.$ 

This identity is the main advantage of the difference setup  
and has no counterpart in the Harish-Chandra theory. The limits of
$p_n^{(k)}$ as $q\to 1$ and $X=e^x$ are spherical functions for
$k=1/2,1,2$ and the Gegenbauer polynomials for arbitrary $k.$
The nonsymmetric polynomials seem to have no direct relations to 
Lie groups and Lie algebras. However the definition is not quite new in 
the representation theory. Actually it is borrowed from affine Hecke 
algebras (but the representation $\l$ is new).

The main references are [O2, M2, C2]. As far as I know, for the first
time the nonsymmetric polynomials appeared in Heckman's 
lectures (the differential
case). 

The first 5 spherical polynomials are:
$$
\eqalignno{
\ep_0 &\ =\ 1,\ \ep_1\ =\  t^{1/2}X,\ \ep_{-1}\ =
\ t^{1/2}({1-tq\over 1-t^2 q}X^{-1}+
{1-t\over 1-t^2 q}X),\cr
&\ep_2\ =\  t ({1-tq\over 1-t^2 q}X^{2}+
{q-qt\over 1-t^2 q}), &\eqnu\cr 
\ep_{-2}&= {t(1-tq)(1-tq^2)\over (1-t^2 q)(1-t^2 q^2)}
 (X^{-2}+
{1-t\over 1-tq^2}X^2+ {(q+1)(1-t)\over 1-tq^2}).
\label\five\eqnum*
}
$$

\vskip 0.2cm
{\bf Duality.}
Let us prove the duality in detail. We will use the anti-involution
$\vph$ from (\ref\phiant) and the following {\it evaluation}
on $ \HH\ :$
$$
\eqalignno{
&\{c_{i\epsilon j}X^i T^\epsilon Y^j\}_t\ =\ 
c_{i\epsilon j}t^{-i/2} t^{\epsilon/2} t^{j/2}.
&\eqnu
\label\evalt\eqnum*
}
$$
It is $\vph$-invariant: $\{\vph(H)\}_t=\{H\}_t.$ Hence
$\{A,B\}_t=\{B,A\}_t$ for the pairing 
$$\{A,B\}_t\equal
\{\vph(A)B\}_t,\ A,B\in \HH.
$$

The evaluation map is  the composition
$$
\eqalignno{
&\HH{\buildrel\alpha\over\longrightarrow}\l\ =\ \C[X,X^{-1}]
{\buildrel\beta\over\longrightarrow}\C,
&\eqnu
\label\evalt\eqnum*
}
$$
where $\alpha$ is a residue modulo the ideal $\HH(T-t^{1/2})+\HH(Y-t^{1/2})$
and $\beta(f(X))=f(X\mapsto t^{-1/2}).$
Given any $f,g\in\C[X,X^{-1}]$, 
$$
\eqalignno{
&\{f(X),g(X)\}_t\ =\ \{\vph(f(X))g(X)\}_t\ =\cr
&\{f(Y^{-1})g(X)\}_t\ =\ f(\widehat{Y}^{-1})(g)(t^{-1}).
&\eqnu
\label\evalsym\eqnum*
}
$$

Here by $\widehat{H}$ we mean the image of $H\in \HH$ in the polynomial
representation $\l$ (with $X$ being $q^x$).
The last relation readily follows from the interpretation of $\l$
as the induced representation from the character 
$\chi(T)=t^{1/2}=\chi(Y)$ of
$\h_Y=\langle T,Y\rangle:$
$$
\eqalignno{
&\l=\hbox{Ind}_{\h_Y}^{\HH}(\chi)
=\HH/\{\HH(T-t^{1/2})+\HH(Y-t^{1/2})\}\simeq\C[X,X^{-1}].
}
$$

Finally, thanks to (\ref\evalsym), (\ref\epspol), 
and the spherical normalization: 
$$
\ep_n(m_\sharp) = \{\ep_n,\ep_m\}_t\ =\  \{\ep_m,\ep_n\}_t=\ep_m(n_\sharp).
$$

The duality can be used to calculate
the renormalization constants, ``$q,t$-dimensions'',
$e_n(-k/2).$ See the formulas
in [C4] (Appendix).

%
%
%
\section{Fourier transform}
Combining $\tau_\pm$  (see (\ref\taupm)), we will prove that
the spherical polynomials multiplied by  the Gaussian are
eigenfunctions of the $q$-Fourier transform and establish 
(\ref\ggmnid). Here the nonsymmetric polynomials are much more
convenient to deal with than the symmetric (even) ones.

We need the unitary structure of the polynomial representation $\l,$
which is given by the
$$
\eqalignno{
&\mu(x;q,k)\ =\ \prod_{j=0}^\infty {(1-q^{j+2x})(1-q^{j+1-2x})\over
(1-q^{j+k+2x})(1-q^{j+k+1-2x})},
&\eqnu
\label\muagain\eqnum*
}
$$
from (7). We extend the map $X^*=X^{-1}, Y^*= Y^{-1},
T^*= T^{-1}, q^*= q^{-1}, t^* t^{-1}$ to an
anti-involution of $\HH.$ Its restriction to  
$\l$ will be denoted by $*$ too.
Recall that $X$ is identified with $q^x$ and $X(z)=q^z.$ 

Since $*$ changes $q,t$ we have to be more precise with the field of
constants. From now on let 
$$
\eqalignno{
&\l\equal\Q_{q,t}[X,X^{-1}] \for \Q_{q,t}\equal\Q(q^{1/4},t^{1/2}).
&\eqnu
\label\poprep\eqnum*
}
$$

Setting $\mu_\circ=\mu/\lan \mu\ran_{\circ},$
$$
\eqalignno{
&\lan H(f), g \ran_\circ\ =\ \lan f, H^*(g)\ran_\circ \for
\lan f, g \ran_\circ\equal \hbox{Const\ Term}( f g^*\mu_\circ ). 
&\eqnu
\label\muform\eqnum*
}
$$
Here $\mu_\circ$ is considered as a Laurent series. Actually its
coefficients are rational functions in terms of $q,t.$ This gives
that $\mu_\circ^*=\mu_\circ$ and $\lan f, g \ran_\circ=
\lan g, f \ran_\circ^*$ for
$f,g\in \l.$ Thus the form is $*$-hermitian.

An immediate application of (\ref\muform) is the  orthogonality
of $e_n$ (and $p_n$) for pairwise distinct indices. The $\mu$-norms 
of the $e$-polynomials and $p$-polynomials are known 
(for arbitrary root systems)
but we do not need them here.  

\vskip 0.2cm
{\bf Functional representation.}
Let us discuss the sharp-variant of the polynomial
representation. The space will be 
$$
\eqalignno{
&\f\equal \hbox{Funct}_{fin}(\ \bowtie\, ,\Q_{q,t}),\ 
\ \bowtie\ =\Z_\sharp=\{n_\sharp, n\in \Z\}.
&\eqnu
\label\funrep\eqnum*
}
$$
By Funct$_{fin}$ we mean finitely supported functions.

The $X$ becomes the operator of  multiplication by $q^{z}$ for
$z\in \ \bowtie\ .$ 
The action of $s,p$ remains the same: $s(f)(z)=f(s(z)),\ $
$p(f)(z)=f(z+1/2).$ So they can't be applied to
any $f\in \f.$ However $sp$ is well-defined. The formula
for $\widehat{T}$ reads 
$$
\eqalignno{
&T'(f)(z)={t^{1/2}q^{2z}-t^{-1/2}\over q^{2z}-1} f(s(z))-
{t^{1/2}-t^{-1/2}\over q^{2z}-1}f(z).
&\eqnu
\label\funti\eqnum*
}
$$
It is well-defined in $\f$  thanks to the factor 
$t^{1/2}q^{2z}-t^{-1/2}$ which vanishes whenever $s(z)\not\in\ \bowtie\ $
for $z\in \ \bowtie\ .$ Therefore $Y=spT$ acts in $\f$ too,
and we have a representation $H\mapsto H'$ of $\HH\ .$  

We set $f^*(z)=f(z)^*,\ $ $\mu_\bullet(z)\equal\mu(z)/\mu(-k/2).$
Then $\mu_\bullet^*=\mu_\bullet.$ 
The inner product is a direct counterpart of $\lan f, g \ran_\circ:$ 
$$
\eqalignno{
&\lan f, g \ran_\bullet\equal \lan fg^* \ran_\bullet \for 
\lan f\ran_\bullet=\sum_{n\in \Z} f(n_\sharp)\mu_\bullet(n_\sharp),
&\eqnu\cr
\label\muonef\eqnum*
&\mu_\bullet(n_\sharp)=\mu_\bullet((1-n)_\sharp)=
q^{-k(n-1)}\prod_{j=1}^{n-1}{1-q^{2k+j}\over 1-q^j} \for n>0.
}
$$
Actually $\lan\ \ran_\bullet$ was already used above in (\ref\sharone).
The old and new ones coincide for even functions $f(x).$ Use the
above formulas for $\mu_\bullet$ to check it.  

The final claim is that $\f$ is $*$-unitary with respect to 
$\lan f,\, g \ran_\bullet.\ $  Cf. (\ref\muform).
The space $\FF$ of all functions is also an $\HH\ $-module.
The pairing can be  
extended to some subspaces of $\FF$ . For instance, 
the scalar products 
$\lan e_n, e_m \ran_\bullet$ are well-defined for $\Re k << 0.$
However only finitely many of them converge for any given $k.$ 
To ``integrate'' all $\ep_n$ for each $k$ we
will add the Gaussians.

\vskip 0.2cm
{\bf Main Theorem.}
Generalizing (\ref\progers), we  calculate the Fourier transforms
of the spherical polynomials times the Gaussian. 
There are several variants of the Fourier transform. We will discuss
only two of them:
$$
\eqalignno{
&F_\circ(f)(n_\sharp)=\lan f(x)\,\ep_n(x)\ran_\circ \and
F_\bullet(f)(n_\sharp)=\lan f(x)\,\ep_n^*(x)\ran_\bullet. 
&\eqnu
 \label\fourtwo\eqnum*
}
$$
Here $f(x)$ is taken from (a proper completion of)
$\l$ and $\f$ respectively,
whereas the Fourier-images belong to $\FF.$

Recall that by $\tilde{\ga}^{-1}$ we mean 
the Gaussian $\ga=q^{x^2}$ considered as a Laurent series
(see (8)). It is assumed that $0<q<1.$
We will use the automorphisms 
$\tau_{\pm}$ from (\ref\taupm).

\proclaim{Theorem} 
i) The transforms $F_\circ$ induces on $\HH\ $ the involution
$\si=\tau_+\tau_-^{-1}\tau_+$ and is unitary on $\l\ :$ 
$$
\eqalignno{
&F_\circ(\widehat{H}(f))\ =\ 
(\si(H))'(F_\circ(f)) \for H\in \HH,\ f\in \l,\cr
&\lan f,g\ran_\circ\ =\ \lan F_\circ(f),F_\circ(g)\ran_\bullet,\ f,g\in \l.
&\eqnu
\label\fourinv\eqnum*
}
$$ 
Respectively, $F_\bullet:\f\to \FF $ induces $\si^{-1}.$ 

ii) Given an arbitrary $k$ and any $n,m\in \Z,$ 
$$
\eqalignno{
&\lan\ep_m\ep_n\tilde{\ga}^{-1}\ran_\circ\ =\
q^{m^2+n^2+2k(|m|+|n|)\over 4}\,\ep_n(m_\sharp)\,
\lan\tilde{\ga}^{-1}\ran_\circ,
&\eqnu\cr
 \label\fgao\eqnum*
&\lan\ep_m\,\ep^*_n\,\ga\ran_\bullet\ =\
q^{-{m^2+n^2+2k(|m|+|n|)\over 4}}\,\ep_n(m_\sharp)\,\lan\ga\ran_\bullet.
&\eqnu
 \label\fgaone\eqnum*
}
$$
\label\FGAUS\theoremnum*
\endproclaim

{\it Proof.} The first part readily results from the
irreducibility of $\l$ and $\f.$
The second part is based on the following fact:
$\widehat{H}\tga^{-1}=\tga^{-1}\tau_+(\widehat{H})$ in the polynomial
representation extended by the Gaussian. Indeed,
the conjugation  by $\tilde{\ga}$
corresponds to $\tau_+$ on the generators 
$X,$ $T,$ and $sp.$
The same  holds in the functional representations $\f$
if $\ga=q^{z^2}$ is treated as an element of $\FF.$
We prove here only (\ref\fgao). 

The map 
$$
F_\ga: f(X)\mapsto \tilde{f}(n_\sharp)=
\ga^{-1}(n_\sharp)\lan f \ep_n\tilde{\ga}^{-1}\ran_\circ,
$$
induces the involution $\tau_+^{-1}\si\tau_+^{-1}=\tau_-^{-1}$
on $\HH\ .$
This map acts from $\l$ to $\FF$ where the field of constants is
extended by  $q^{k^2/4}.$ 

The automorphism $\tau_-$ fixes  $Y$. Hence
the image of $f=\ep_n$ is an eigenfunction of the discretizations
$Y'$ with the same eigenvalue. Let us prove that 
$F_\ga(\ep_n)$ has to be proportional
to $\ep'_n=\ep_n(m_\sharp).$ 

One may assume that $n=0$ (an exercise).  
The function $g=F_\ga(1)$ satisfies
$t^{-1/2}T'(g)=$ $g=(sp)(g).$ 
We have used here that $\tau_-$ fixes $T$ and $sp.$ 
Thus $g$ is invariant with respect to $s$ and $p,$ 
which means that it is a constant.  

Setting $F_\ga(\ep_n)=
h_n\ga(n_\sharp)\,q^{-k^2/4}\,\ep'_n$ for $h_n\in \Q_{q,t},$
we need to check that $h_n=1.$ It is true of course for $h_0.$ 
However
$h_m=h_n$ for all $m,n\in \Z$
because both the left and right sides of (\ref\fgao) are 
$m\leftrightarrow n$ symmetric. As
to the right-hand side, it is because of the duality.
\proofbox

It follows from the theorem, that 
the composition $F_\bullet F_\circ$ is nothing else but the
discretization map $f\mapsto f'=f(n_\sharp).$ Indeed, 
$\F_\circ$ is an isomorphism $\l\to \f$ and
$F_\bullet F_\circ$ induces the identity on $\HH\ .$ 
So it has to be proportional
to $f\mapsto f'$ because of the irreducibility of $\l.$
It suffices to calculate the coefficient of 
proportionality for $f=\ep_0=1,$ which is  $1$ thanks to the normalization
of $\mu_\circ$ and $\mu_\bullet.$  

Explicitly, $F_\circ$ sends 
$$
\eqalignno{
&\ep^*_m\mapsto \mu_\bullet^{-1}(m_\sharp)\de_{m_\sharp}\for
\de_{m_\sharp}(n_\sharp)=\de_{m_\sharp n_\sharp}.
&\eqnu
 \label\fouriso\eqnum*
}
$$  
The simplest way to see it is to apply $F_\bullet$ to 
$\de_{m_\sharp}$ and to use that $F_\bullet F_\circ=$id.

This observation is directly related to the Plancherel
formula for the $p$-adic spherical transform. The nonspherical case
requires a variation of the initial point in the Jackson integral.
It is  $-k/2$ in this paper. See [C1],[C2]. 
Analytic problems seem not very difficult, at least
in the one-dimensional case [KS]. 
It is likely that there are connections
with $p$-adic Lusztig's theory and recent results due to
Heckman, Opdam on the general Plancherel formula for affine 
Hecke algebras. Hopefully the double Hecke algebra can give here
more than just a $q$-deformation. 

The analytic considerations become much simpler if we
switch from $\l$ to $\tga^{-1}\l$ and adjust properly the functional
case. However the Gaussian  collapses in the $p$-adic limit:
$q\to \infty$ with $X,t$ being fixed. The zonal limit sending  $q\to 1$
and fixing  $X,k$ destroys the Gaussian too. 
If the Gaussian partially survives under such limits then it may lead to
a new technique in the classical harmonic analysis,
but there are no confirmations so far.

%
%
%
\section{Roots of unity}
Let $q=\exp(2\pi i/N)$.  Actually the formulas which do
not contain the imaginary unit $i$
explicitly  hold for any primitive
$q.$ Indeed, we can apply the Galois automorphisms. 
However we need ``the least'' $q$  to ensure the positivity
of the inner product. The sign of $q^{1/2}=\pm \exp(\pi i/N)$
is sometimes important too. As to  $q^{1/4},$ the formulas are true 
for either choice. We remind the reader that by $q^{n/4}$ we mean
$(q^{1/4})^n$ unless otherwise stated. 

We begin with the {\it main sector}  which is 
$0<k< N/2,\ k\in \Z.$
Recall that $q^x(m/2)=q^{m/2}$ for 
$m\in \Z$. The field of constants is $\Q_q=\Q(q^{1/4}).$ 
The involution $*$ becomes the complex conjugation
$(q^{1/4})^*=q^{-1/4},$ trivial on $x.$
For $ -N< m\le N,$  
$$
\eqalignno{
&\bowtie'\ \equal \{m_\sharp\mid \mu_\bullet(m_\sharp)\neq 0\}\ =\cr
&\{{-N+k+1\over 2},\ldots, -{k\over 2},\ {k+1\over 2}, 
\ldots,{N-k\over 2}\}.
&\eqnu
 \label\bowtiep\eqnum*
}
$$
We use (\ref\muonef):
$$
\eqalignno{
&\mu_\bullet(m_\sharp)=\mu_\bullet((1-m)_\sharp)=
q^{-k(m-1)}\prod_{j=1}^{m-1}{1-q^{2k+j}\over 1-q^j} \for m>0.
}
$$


The space 
$\f'=\f'(k)=\hbox{Funct}(\ \bowtie'\,,\Q_q)$ has a unique structure of
a $\HH$-module making the {\it discretization map}
$$\l \ni f\ \mapsto\ f'=f(m/2)\in \f'
$$ 
a  $\HH$-homomorphism. Its dimension
is $2(N-2k).$
Indeed, we need to check that $T'$ from (\ref\funti)
and $sp$ are well-defined on $\f',$ which follows directly from
the definition of $\ \bowtie'.$  
The conceptual proof is as follows.

Setting
$$
\lan f\ran'\equal 
\sum_{m=-N+1}^{N} f(m_\sharp)\mu_\bullet'(m_\sharp),
$$
the form $\lan f,g\ran'=$ 
 $\lan fg^*\ran'$
is $*$-hermitian
on $\l.$ The proof is the same as
for $\lan f,g\ran_\circ.$ We utilize the defining difference equation
of $\mu.$ 
The module $\f'$ is nothing
else but $\l$ modulo the radical of this form. Here we need to
exclude the {\it special case} of odd $N$ and  $q^{1/2}=-\exp(\pi i/N),$
when the image of $\p$ is two times smaller than  $\f'.$

The same hermitian form when considered on
$\f'_k$ make it is a $*$-unitary representation. 
Moreover the inner product is positive thanks to the special
choice $q=\exp(2\pi i/N).$
Here the sign of  $q^{1/2}=\pm\exp(\pi i/N)$ can be arbitrary.
The inner product does not involve $q^{1/2}.$ So 
the special case is included.

The module  $\f'$ is irreducible
if $N$ is even or  if $q^{1/2}=\exp(\pi i/N)$ for odd $N.$
Indeed, $1$ is a unique eigenvector of
$Y'$ with the eigenvalue $t^{1/2}=q^{k/2}.$ On the other hand,
it is cyclic, i.e. it generates $\f'$ as a $\HH$-module. 
This and the semisimplicity of $Y'$ result in the irreducibility.
In the special case $N=2M+1,\ $ $q^{1/2}=-\exp(\pi i/N),$ there are two 
irreducible components:
$$
\eqalignno{
&\f'\ =\ (\f')^0\oplus\f^1,\ \
\hbox{dim} (\f')^j=N-2k,
&\eqnu\cr 
\label\evenni\eqnum* 
&(\f')^j\ =\
\{f\in\f' \ \mid \ f(n_\sharp+N/2)=(-1)^{j}f(n_\sharp),\ n_sharp<0 \}.
}
$$
The component $(\f')^0$ is the image of $\p.$ 

It is not difficult to calculate the spectrum of $Y'$ exactly
when $\f'$ is irreducible.
The images $\ep'_n$ of the spherical polynomials $\ep_n$
are well-defined when $n_\sharp\in \ \bowtie'\, .$ They are 
$Y$-eigenfunctions, namely
$Y'(\ep'_n)=q^{-n_\sharp}\ep'_n,$ and linearly generate $\f'.$  
We have the duality: $\ep'_m(n_\sharp)= \ep'_n(m_\sharp).$

In the case of (\ref\evenni), we get the $Y$-spectrum in  $(\f')^0$
taking even $m,n.$ We will not discuss the diagonalization of $Y$
in  $(\f')^1$ here.

The $\f'$ inherits all properties of $\l,\f.$ 
For instance,  $PSL^c_2(\Z)$ acts on it.
The automorphism $\tau_+$
is simply the multiplication by the image $\ga'\in \f'$
of the Gaussian $\ga=q^{x^2}.$

From now on we skip the prime in the formulas: $\ep_n$ and $\ga$ will
be always considered in the corresponding functional
representation. 

Let $\f'=\f'(k)=(\f')^+(k)\oplus (\f')^-(k)$ where
$T=\pm t^{\pm 1/2}$ on $(\f')^{\pm}.$
The corresponding dimensions are
dim$\f'=2(N-2k)=$ $ (N-2k+1)+(N-2k-1)$ provided that $k<N/2.$  
The components $(\f')^\pm$ are $PSL^c_2(\Z)$-invariant.
Calculating the  action of  $PSL^c_2(\Z)$
in  $(\f')^+$  we come to the formulas
from [Ki,C4] as $q^{1/2}=\exp(\pi i/N).$  
The $PSL^c_2(\Z)$-module  $(\f')^-(k)$ is 
$PSL^c_2(\Z)$-isomorphic to  $(\f')^+(k+1)$.
The $(\f')^+(1)$ is the {\it Verlinde algebra}. 
It is a subalgebra of $\f'$ (but not a submodule).

\vskip 0.2cm
{\bf Gaussian sums.}
Let us now adjust the Main Theorem to $\f'.$
Both variants of the Fourier transform considered above may be used:
$$
\eqalignno{
&F_\circ(f)(n_\sharp)=\lan f \ep_n\ran' \and
F_\bullet(f)(n_\sharp)=\lan f\,\ep_n^*\ran'. 
&\eqnu
 \label\fourtwor\eqnum*
}
$$
The summation in $\lan\ \ran'$ 
is  over $\ \bowtie'.$ 
Here we assume that $f\in (\f')^0$ in the special case of (\ref\evenni).
Otherwise the transforms are zero.

They induce the involutions $\si$ and $\si^{-1}$ on $\HH\ $
and are unitary with respect to $\lan f,g \ran'.$
Here $f,g\in \f'.$ 

\proclaim{Theorem} 
Given  $n,m$ such that $n_\sharp, m_\sharp\in \ \bowtie'\, ,$ 
$$
\eqalignno{
&\lan\ep_m\ep_n\,\ga^{-1}\ran'\ =\
q^{m^2+n^2+2k(|m|+|n|)\over 4}\, \ep_n(m_\sharp)\, (C'_k)^*,
&\eqnu\cr
 \label\fgaor\eqnum*
&\lan\ep_m \,\ep_n^*\,\ga'\ran'\ =\
q^{-{m^2+n^2+2k(|m|+|n|)\over 4}}\, \ep_n(m_\sharp)\, C'_k.
&\eqnu
 \label\fgaonr\eqnum*
}
$$
Following (\ref\sharpga):
$$
\eqalignno{
&C'_k\ =\ \lan \ga\ran'\ =\ 
\prod_{j=1}^k(1-q^j)^{-1}
\sum_{m=-N+1}^{N}q^{m^2/ 4}.
&\eqnu
\label\ctgar\eqnum*
}
$$
\label\FGAUSR\theoremnum*
\endproclaim

The last formula coincides with (10) for $q^{1/2}=\exp(\pi i/N)$. 
It can be $0=0$ as $N=2M+1,q^{1/2}=-\exp(\pi i/N).$
This happens exactly when $N=4L+1,$ because the Gaussian sits
in $(\f')^1$ for such $N$
and is orthogonal to $1$ and all spherical polynomials.
For  $N=4L+3,$ it belongs to  $(\f')^0$ (see (\ref\evenni). 

We note that (\ref\ctgar) can be
deduced from (\ref\sharpga) using the standard limiting procedure
from the $\theta$-series to the Gaussian sums.

There are several other classical formulas similar to  (10).
To cover them all we need to diminish the double Hecke algebra
and its irreducible representations.

\vskip 0.2cm
{\bf Generalized Gaussian sums.}
Let $\widetilde{\HH}$ be a subalgebra of $\HH\ $ generated by
$X^2,T$ and $Y^2.$ We will assume for a while that $q$ and $t$
are generic. The image of $T_0\equal Y^2T^{-1}$ in the polynomial 
representation $\l$ can be readily calculated:
$$
\eqalignno{
&\widehat{T_0}=
t^{1/2}s_0+(q^{1-2x}-1)^{-1}(t^{1/2}-t^{-1/2})(s_0-1), \ s_0=sp^2.
&\eqnu
\label\toop\eqnum*
}
$$
This gives that  $\widetilde{\HH}$ satisfies the PBW theorem, which
can be of course checked directly using the abstract relations between
the generators. All the symmetries of $\HH\ $ hold for  $\widetilde{\HH}.$
Moreover, we can replace $\l$ by its ``even''
part $\l^0\equal \Q_{q,t}[X^2,X^{-2}],$ which is an irreducible
$\widetilde{\HH}$-module. The odd part $\l^1\equal X\l^0$ is an irreducible
module as well. Respectively we may consider the spherical polynomials
$\ep_n$ either for even or odd $n.$ They  
generate  $\l^0$ and $\l^1$ respectively.

Similarly, we diminish $\f:$ 
$\ \widetilde{\f}\equal \hbox{Funct}_{fin}((2\Z)_\sharp,\Q_{q,t}).$
Also the summation in the definition of $\lan \ \ran_\bullet$ will be
over the set $(2\Z)_\sharp.$  We write $\lan \ \ran_{\bullet\!\bullet}.$

Formula (\ref\fgaone) holds for  $\lan \ \ran_{\bullet\!\bullet}$ 
with arbitrary $m,n.$ 
The identity (\ref\sharpga) now reads:
$$
\eqalignno{
&\lan q^{x^2}\ran_{\bullet\!\bullet} =
q^{k^2/4}\sum_{j=0}^\infty {1-q^{2j+k}\over 1-q^{k}} \prod_{l=1}^{2j}
{1-q^{l+2k-1}\over  1-q^{l}} q^{j^2-jk}=\cr
&q^{k^2/4}\prod_{j=1}^\infty
{1-q^{j+k}\over 1-q^j}\sum_{j=-\infty}^{\infty}q^{j^2-jk}.
&\eqnu
\label\sharpgat\eqnum*
}
$$

We turn to the roots of unity:
$q=\exp(2\pi i/N)$ where $0< k\le M=[N/2],$
$q^{1/2}=\pm\exp(\pi i/N).$ 
The sign of $q^{1/2}$ can be arbitrary. All formulas hold for either
choice.

We get the identity
$$
\eqalignno{
&\sum_{j=0}^{M-k} {1-q^{2j+k}\over 1-q^{k}} \prod_{l=1}^{2j}
{1-q^{l+2k-1}\over  1-q^{l}} q^{j^2-jk}=\cr
&\prod_{j=1}^k
{1\over (1-q^j)}\sum_{j=0}^{N-1}q^{j^2-jk}.
&\eqnu
\label\shgatr\eqnum*
}
$$ 

This identity may be $0=0$. It happens when $N=2M$ and 
$M-k$ is odd. Setting $k=M,$   
$$
\eqalignno{
&\sum_{j=0}^{N-1}q^{j^2} = q^{L^2}\prod_{j=1}^M (1-q^j)
\for N=2M+1, L=M/2\mod N,\cr
&\sum_{j=0}^{M-1}q^{j^2}(-1)^j\ =\ \prod_{j=1}^{M-1} (1-q^j) \for N=2M.
&\eqnu
\label\gatr\eqnum*
}
$$ 
The latter gives a product formula for
the so-called {\it generalized Gaussian sum}. When
considered for all primitive roots $q,$ it 
plays the key role in the quadratic reciprocity. 

Using that $\sum_{j=0}^{N-1} q^{(j-k/2)^2}$ does not depend on
$k,$ we may simplify (\ref\shgatr)  under the assumption
that $M-k$ is even for $N=2M\ :$
$$
\eqalignno{
&\sum_{j=0}^{M-k} {1-q^{2j+k}\over 1-q^{k}} \prod_{l=1}^{2j}
{1-q^{l+2k-1}\over  1-q^{l}} q^{j^2-jk}=\cr
&q^{(L-K)(L+K)}\prod_{j=k+1}^M
(1-q^j)\for L\pm K= (M\pm k)/2 \mod N.
&\eqnu
\label\shoe\eqnum*
}
$$ 
When $M-k$ and $N$ are even the left-hand side is zero.

Let us establish the counterpart of Theorem \ref\FGAUSR and
as a by-product clarify why the case of even $N$ and odd $M-k$
is exceptional. We assume that $k<N/2.$ Now
$\ \bowtie''\ \equal$
$$
\eqalignno{
&\{{-N+k\over 2}+1,\ldots, -{k\over 2},\ {k\over 2}+1, 
\ldots,{N-k\over 2}\} \if N=2M,
&\eqnu\cr
 \label\bowtiee\eqnum*
&\{{-N+k+1\over 2},\ldots, -{k\over 2},\ {k\over 2}+1, 
\ldots,{N-k-1\over 2}\} \if N=2M+1.
}
$$

So the dimension of the  $\widetilde{\HH}$-module
$\f''\equal\hbox{Funct}(\ \bowtie''\,,\Q_q)$ 
is $N-2k.$ It is irreducible for odd $N$ and has
two irreducible $\widetilde{\HH}$-components:
$$
\eqalignno{
&\f''\ =\ (\f'')^0\oplus(\f'')^1,\ \
\hbox{dim} (\f'')^j=M-k \for N=2M:\cr
&(\f'')^j=\hbox{Image}(\l^j)=
\{f\in\f'' \ \mid \ f(n_\sharp+M)=(-1)^{j}f(n_\sharp)\},
&\eqnu
 \label\evenn\eqnum*
}
$$
where we take  negative $n_\sharp.$ 

The Gaussian $q^{n_\sharp^2}$ belongs to
$(\f'')^0$ precisely when $M-k$ is even. Otherwise it
sits in $(\f'')^1.$ This explains the degeneration of
(\ref\shgatr). Note that $\f''$ for odd $N$ is nothing else but
$(\f')^0$ for $q^{1/2}=-\exp(\pi i/N).$ So $(\f')^0$ from (\ref\evenni)
remains irreducible upon the restriction to $\widetilde{\HH}\ .$

Let us adjust (\ref\fgaor). We need
$$
\eqalignno{
&\lan f\ran''\equal \sum f(z)\mu_\bullet(z),\ z\in  
\ \bowtie'', \and \cr 
&C''_k\equal 
q^{L^2}\prod_{j=k+1}^M
(1-q^j), \where
&\eqnu\cr
\label\shiodo\eqnum*
& L= M/2\ (N=2M),\ L=M/2\mod N\ (N=2M+1). 
}
$$ 
The latter means that $q^{L^2}=i^{-M^2N}q^{M^2/4}$ in the
case $N=2M+1$ if $q^{1/4}=+\exp(\pi i/(2N)).$
The plus-sign in the latter is just the matter of normalization.

We pick  $1\le n,m\le N-2k.$ 
The images of $\ep_n$ in $\f''$
are linearly independent in this range and linearly
generate the whole space. 

Finally for the same choice $q^{1/4}=\exp(\pi i/(2N))$ as above:
$$
\eqalignno{
&\lan\ep_m\, \ep_n^*\,\ga\ran''\ =\
C\, q^{-{m^2+n^2+2k(|m|+|n|)\over 4}}\, \ep_n(m_\sharp)\, C''_k,
&\eqnu\cr
 \label\fgaonrr\eqnum*
&C=i^{N(m^2+n^2+2k(|m|+|n|))}\if N=2M+1,\cr 
&C=0 \if M+k+n+m \hbox{\ is\ odd} \for N=2M,\cr
&C\ =\ 1\ \hbox{\ otherwise\ .}
}
$$

\vfil
It gives  an explicit description
of the Fourier transform 
$F_\bullet$ (see (\ref\fourtwor)).
Here and further we will skip the analogous
formulas for $F_\circ.$ Cf. (\ref\fgaor).
The following changes are
required to go from $F_\bullet$ to $F_\circ:$

\item{a)} $\ga\mapsto \ga^{-1},\ \ep_n^* \mapsto \ep_n\ $ in $\ \lan\ \ran$, 
\item{b)} $q^{-(m^2+n^2+2k(|m|+|n|))/ 4}\, \mapsto\,
q^{(m^2+n^2+2k(|m|+|n|))/ 4}\ ,$
\item{c)} the constants $C$ are changed by their conjugations $C^*.$

%
%
%
\section{Half-integral k}
Let us discuss now the case which is different from the
classical theory of Gaussian sums. Actually the above considerations
are already beyond the classical theory thanks to $k.$ However
when $k=1$ and  for other integral $k>1$ there are strong links to 
the classical formulas, as we demonstrated above. We are going  to
consider half-integral $k,$ where the classical origin is less
clear. For instance, the Gaussian sum for $k=1/2$ is 
of collapsing type, in contrast to  $k=1.$ Recall that $k=1/2$ is
the case of $SL_2(\R)/SO_2$ in the 
the Harish-Chandra theory.

This does not mean  that the formulas below, especially
those without  spherical polynomials, cannot be verified by elementary 
methods. It is always probable in the one-dimensional case. 
Moreover, half-integral negative $k$ are directly connected
with integral positive $k$ for odd $N,$ as we will see.  

We note that the formulas of this sections combined with
the previous ones give a complete list of Fourier
transforms of self-dual spherical irreducible
unitary representations of double affine Hecke algebras at roots
of unity (see [C5]).  

Let $k\in 1/2+\Z.$ We follow the standard procedure and
begin with generic $q.$
Using (9),
$$
\eqalignno{
&\sum_{j=0}^\infty q^{{(j-k)^2\over 4}}
{1-q^{j+k}\over 1-q^{k}} \prod_{l=1}^j
{1-q^{l+2k-1}\over
 1-q^{l}}\cr
&=\  2q^{1/16} \prod_{j=1}^\infty
{1-q^{j+k}\over 1-q^{-1/2+j}}
&\eqnu\cr
\label\shhalf\eqnum*
&=\  2q^{1/16}\prod_{j=0}^{s}
{1\over 1-q^{1/2+j}}\for k=1/2+s,\ s\in \Z_+,\cr
&=\  2q^{1/16}\prod_{j=1}^{s}
(1-q^{1/2-j})\for k=-1/2-s,\ s\in \Z_+.\cr
}
$$

Now let $q=\exp(\pi i/N).$ The 
definition of $\ \bowtie'\ $ and
the analysis of the representations of $\HH\ $ for
$0<k<N/2$ are practically the same. However
we need to assume that $q^{1/2}=\exp(\pi i/2)$ from the
very beginning to provide
the positivity of the inner product. Recall that 
the values of $\mu_\bullet$ on $\ \bowtie'\ $ 
did not involve $q^{1/2}$ for integral $k.$

The definition of the
$\HH$-module $\f'$ remains unchanged. It is irreducible
of dimension $2(N-2k).$ The summation
in (\ref\shhalf) is from $j=0$ to $j=N-2k.$
There is only one change in Theorem \ref\FGAUSR. 
We need  to replace $C_k'$ by 
$$
\eqalignno{
&C_k'\ =\   2q^{1/16}\prod_{j=0}^{k-1/2}
{1\over (1-q^{1/2+j})}.
&\eqnu
\label\shhalfc\eqnum*}
$$

\vskip 0.2cm
{\bf The negative case.}
One may also consider the interval $N\le 2k< 2N$
for half-integral $k.$ For such $k,$ 
the space of functions on 
$$
\eqalignno{
&\{{k+1\over 2},{k+2\over 2},\ldots, N-{k\over 2}\}
&\eqnu
 \label\bowtil\eqnum*
}
$$
is an irreducible
$\HH$-module of dimension $2N-2k.$ 

It is more convenient
to switch from  $k$ to $\bar{k}=k-N.$ 
This may change the sign of $t^{1/2}$ in the definition
of the double Hecke algebra but the impact is easy to control.
The map $T\mapsto -T,\ Y\mapsto -Y$
fixing $ X,q,t$ 
is an automorphism of $\HH\ .$ 

Thus let us take $\bar{k}=-1/2-s$ instead of $k$ and set 
$$
\eqalignno{
&\bar{\bowtie}' \equal 
\{{\bar{k}+1\over 2},\ldots,-{\bar{k}+1\over 2},-{\bar{k}\over 2}\}=
&\eqnu\cr
 \label\bowtill\eqnum*
&\{{\bar{k}+1\over 2},\ldots,-{1\over 4},{1\over 4}=
{\bar{k}+s+1\over 2},\ldots,
-{\bar{k}\over 2}={\bar{k}+2s+1\over 2}\}.
}
$$
Then $\bar{\f}'= \hbox{Funct}(\ \bar{\bowtie}'\, ,\, \Q_q)$
is an irreducible $\HH(\bar{k})$-module of the same dimension
 $-2\bar{k}=2N-2k.$ We put $\HH(k)$ instead of $\HH\ $ to show
explicitly the dependence on $k.$

Let $q=\exp(2\pi i/N)$ and, moreover, $q^{1/2}=-\exp(\pi i/N).$
The minus sign is necessary in the case under consideration
to make the inner product positive. All  formulas will
hold for either sign, as well as the analysis of representations. 

The function   $\mu_\bullet$ 
(it is the same for $k$ and $\bar{k}$) makes
 $\bar{\f}'$  unitary. Recall that
$\lan\ f,g \ran'$ is the summation of $fg^*\mu_\bullet$
over the set $\ \bar{\bowtie}'\ .$ This form is hermitian and positive.

The module  $\bar{\f}'$ is  generated by the eigenfunctions
of $Y',$ which are the images of $\ep_n$ for $n_\sharp\in 
\ \bar{\bowtie}'\ .$ Generally speaking, the justification
of the  existence of  $\ep_n$
requires either the technique of intertwiners or the Pieri rules.
In the one-dimensional case, the Appendix of [C4] is sufficient.
Then the linear independence of their images in  $\bar{\f}'$ results
from the consideration of the $Y$-eigenvalues. We will skip the
detail (here and further). 
 
Finally,
$$
\eqalignno{
&\lan\ep_m\, \ep^*_n\, \ga\ran'\ =\
q^{-{m^2+n^2+2\bar{k}(|m|+|n|)\over 4}}\, \ep_n(m_\sharp)\, C'_{\bar{k}},
\ \ C'_{\bar{k}}\ =
&\eqnu\cr
 \label\fgaonrhn\eqnum*
&\sum_{j=s+1}^{2s+1} q^{{(j-\bar{k})^2\over 4}}
{1-q^{j+\bar{k}}\over 1-q^{\bar{k}}} \prod_{l=1}^j
{1-q^{l+2\bar{k}-1}\over  1-q^{l}}=  
q^{1/16}\prod_{j=1}^{s}
(1-q^{1/2-j}).
}
$$

\vskip 0.2cm
{\bf The case of $\widetilde{\HH}$\ .} 
As always, we begin with generic identities:
$$
\eqalignno{
&\sum_{j=0}^\infty q^{j^2-kj}
{1-q^{2j+k}\over
1-q^{k}} \prod_{l=1}^{2j}
{1-q^{l+2k-1}\over
 1-q^{l}}\cr
&=\ 
\prod_{j=1}^\infty
{(1+q^{-k-1+2j})\over
(1+q^{k+2j})}
\prod_{j=1}^\infty
{(1-q^{2k+2j})\over
(1-q^{-1+2j})}
&\eqnu\cr
\label\shhalft\eqnum*
&=\  \prod_{j=0}^{s}
{1+q^{k-2j}\over 1-q^{2k-2j}}\for k=1/2+s,\ s\in \Z_+,\cr
&=\  \prod_{j=1}^{s}
{1-q^{2k+2j}\over 1+q^{k+2j}}\for k=-1/2-s,\ s\in \Z_+.\cr
}
$$

Now we switch to  $q^{1/2}=\exp(\pi i/N),$
provided that  $0<2k<N.$
The analysis of $\f''$ in this case is similar to 
that for integral $k.$ However the roles of even and
odd $N$ are inverse.
We set
$\ \bowtie''\ \equal$
$$
\eqalignno{
&\{{-N+k\over 2}+1,\ldots, -{k\over 2},\ {k\over 2}+1, 
\ldots,{N-k\over 2}\} \if N=2M+1,
&\eqnu\cr
 \label\bowtieeh\eqnum*
&\{{-N+k+1\over 2},\ldots, -{k\over 2},\ {k\over 2}+1, 
\ldots,{N-k-1\over 2}\} \if N=2M.
}
$$

The dimension of the  $\widetilde{\HH}$-module
$\f''=\hbox{Funct}(\ \bowtie''\,,\Q_q)$ 
remains  $N-2k$ as in the integral case.
It is irreducible for even $N$ and has the 
$\widetilde{\HH}$-decomposition
$$
\eqalignno{
&\f''\ =\ (\f'')^0\oplus(\f'')^1,\ \
\hbox{dim} (\f'')^j=M-k \for N=2M+1:\cr
&(\f'')^j=\hbox{Image}(\l^j)=
\{f\in\f'' \ \mid \ f(n_\sharp+N/2)=(-1)^{j}f(n_\sharp)\},
&\eqnu
 \label\evenna\eqnum*
}
$$
where $n_\sharp<0.$  
The reduction of (\ref\shhalft) reads:
$$
\eqalignno{
&\sum_{j=0}^{[N/2-k]} q^{j^2-kj}
{1-q^{2j+k}\over
1-q^{k}} \prod_{l=1}^{2j}
{1-q^{l+2k-1}\over
 1-q^{l}}\cr
&=(1+i^{2k-N})\,C''_k,\, \ C''_k\equal \prod_{j=0}^{k-1/2}
{1+q^{k-2j}\over 1-q^{2k-2j}}.
&\eqnu
\label\shhalftr\eqnum*
}
$$
The main formula becomes:
$$
\eqalignno{
&\lan\ep_m\, \ep^*_n\,\ga\ran'' \equal
\sum_{z\in \bowtie''} 
\mu_\bullet(z)\ep_m(z)\, \ep^*_n(z) \ga(z)\ =
&\eqnu\cr
 \label\fgaonrhep\eqnum*
&(1+i^{2k+|m|+|n|-N})\, C''_k\,
q^{-{m^2+n^2+2k(|m|+|n|)\over 4}}\, \ep_n(m_\sharp).
}
$$

\vskip 0.2cm
{\bf The $\widetilde{\HH}\ $-negative case.}
Let us consider  $\bar{k}=k-N$ provided  $N\le 2k<2N.$
I.e.  $-N/2\le \bar{k}=-1/2-s <0$
where $s\in \Z_+$ and  $0\le s<N/2.$

We need to take $q^{1/2}=-\ep(\pi i/N)$ to ensure the positivity
of the hermitian form. However all  formulas below will hold for either
sign of  $q^{1/2}.$

The reduction of  (\ref\shhalft) is straightforward:
$$
\eqalignno{
&\sum_{j=0}^{s} q^{j^2-\bar{k}j}
{1-q^{2j+\bar{k}}\over
1-q^{\bar{k}}} \prod_{l=1}^{2j}
{1-q^{l+2\bar{k}-1}\over
 1-q^{l}}\cr
&= \prod_{j=1}^{s}
{1-q^{2\bar{k}+2j}\over 1+q^{\bar{k}+2j}}.
&\eqnu
\label\shhalftrn\eqnum*
}
$$

Actually this formula is not new. It coincides
with the $C$-part of (\ref\fgaonrhn):
$$
\eqalignno{
&\sum_{j=s+1}^{2s+1} q^{{j^2-2\bar{k}j\over 4}}
{1-q^{j+\bar{k}}\over 1-q^{\bar{k}}} \prod_{l=1}^j
{1-q^{l+2\bar{k}-1}\over  1-q^{l}}=  
q^{1/16-\bar{k}^2/4}\prod_{j=1}^{s}
(1-q^{1/2-j}).
}
$$
The coincidence of the right-hand sides is an elementary exercise:
$$
\eqalignno{
&q^{1/16-\bar{k}^2/4}\prod_{j=1}^{s}
(1-q^{1/2-j})=\prod_{j=1}^{s}
{1-q^{2\bar{k}+2j}\over 1+q^{\bar{k}+2j}}.
&\eqnu
 \label\bowiden\eqnum*
}
$$
Here $q$ is of course arbitrary. This formula is of a certain 
importance in the theory of $\eta$-like identities (as $s\to \infty$).

The explanation is simple.
The $\HH(\bar{k})\ $-module $\bar{\f}'$ introduced above
for negative $\bar{k}=k-N$ 
remains irreducible upon the restriction to $\widetilde{\HH}(\bar{k}).$

The standard $\bowtie''$-set constructed for $\widetilde{\HH}(\bar{k})$
and responsible for the structure of formula (\ref\shhalftrn) is 
$$
\eqalignno{
&\bar{\bowtie}'' \equal 
\{-{\bar{k}\over 2}-s,\ldots,\, -{\bar{k}\over 2},\,
 {\bar{k}\over 2}+1,\ldots,\,{\bar{k}\over 2}+s
\}.
&\eqnu
 \label\bowtilln\eqnum*
}
$$
Here the points from the  second half ($\bar{k}/2+1,$ etc.) sit
between consecutive
pairs of points from the first half.
It is nothing else but a rearrangement of $\bar{\bowtie}'\ $
from (\ref\bowtill). 

The main formula is equivalent to (\ref\fgaonrhn) and  reads 
as follows:
$$
\eqalignno{
&\sum_{j=1}^{s} (q^{-(2j-1)\bar{k}} \prod_{l=1}^{2j-1}
{1-q^{l+2\bar{k}}\over  1-q^{l}}) \,
q^{j^2+\bar{k}j}\,
\ep_m({2j+\bar{k}\over 2})\,\ep^*_n({2j+\bar{k}\over 2})\ +\cr
&\sum_{j=0}^{s} (q^{-2j\bar{k}} \prod_{l=1}^{2j}
{1-q^{l+2\bar{k}}\over  1-q^{l}})\, 
q^{j^2+\bar{k}j} \,
\ep_m(-{2j+\bar{k}\over 2})\,\ep^*_n(-{2j+\bar{k}\over 2})\ =\cr
&q^{-{m^2+n^2+2\bar{k}(|m|+|n|)\over 4}}\, \ep_n(m_\sharp)
 \prod_{j=1}^{s}
{1-q^{2\bar{k}+2j}\over 1+q^{\bar{k}+2j}}.
&\eqnu
\label\strnpo\eqnum*
}
$$
Here  $m_\sharp, n_\sharp\in 
\ \bar{\bowtie}''\ .$

\vskip 0.5cm
{\bf Deformations of Verlinde algebras.}
We will conclude these notes with the following observation.
When $N=2M+1,$ the case $\bar{k}=-1/2-s$ is equivalent
to the case of integral $\hat{k}\equal M-s.$ Here we may replace $-1/2$ 
by $M$ modulo $N$ because $q^{1/2}=-\exp(\pi i/N)$ is a $N$-th root of
unity. To be more exact, $\widetilde{\HH}(\hat{k})$ coincides with 
$\widetilde{\HH}(\bar{k})$ for such $q^{1/2},$ since
the quadratic relation remains the same,  
and  $\f''(\hat{k})$
(see (\ref\bowtiee)) is isomorphic to 
$\f''(\bar{k}).$
The dimensions are of course the
same: dim$\f''(\bar{k}) = -2\bar{k} = 2s+1=$ 
$N-2\hat{k}=$ dim$\f''(\hat{k}).$ 

So we hit the main sector $0<\hat{k}<N/2, \ \hat{k}\in \Z.$
Recall  that the irreducibility of  $\f''(\hat{k}),$
the positivity of the inner product, and all formulas were independent
of the sign of $q^{1/2}$ for such $k.$

Comparing (\ref\shhalftrn) and (\ref\shoe)
we arrive at the identity 
$$
\eqalignno{
&\prod_{j=1}^{s}
{1-q^{1-2j}\over 1+q^{M-s+2j}}\ =\
q^{-s(s+1)/4}\prod_{j=0}^{s-1}(1-q^{M-j}),
&\eqnu
\label\shhide\eqnum*
}
$$
which readily follows from (\ref\bowiden).

It is important that  all claims and formulas about
$\f''(\bar{k}) $ hold for any $q.$
For instance, its dimension is always $2s+1=-2\bar{k},$
the formula (\ref\strnpo) works, and so on. 
The particular choice of
$q$ does not matter. The constraint $k=-1/2-s$ is sufficient.
For $|q|\neq 1,$ the conjugation
$(q^{1/4})^*=q^{-1/4}$ must be  understood formally.

If one  wishes to make the $\mu$-form hermitian and 
positive it is necessary to be specific.
The conditions $q^{1/2}=-\exp(\pi i \om)$ and $ 0<\om s< 1/2$
are sufficient. However
it is not necessary to assume that $q$ is a root of unity.    

One can use  $\HH\ $ as well to establish the 
$\bar{k}\leftrightarrow\hat{k}$ correspondence for odd $N$ and
$q^{1/2}=-\exp(\pi i/N).$
Namely, $(\f')^0(\hat{k})$ is equivalent to  $(\f')(\bar{k}).$
We will remind the reader, that both modules
remain irreducible and become respectively
$\f''(\hat{k})$ and $\f''(\bar{k})$ upon the restriction to 
$\widetilde{\HH}\ .$ Thus the usage of  $\HH\ $ instead
of  $\widetilde{\HH}\ $ does not add anything new.    

We see that the theory of $\f''(k)$ for integral $0<k<N/2$
allows a $q$-deformation when $N=2M+1.$ This covers the {\it little}
Verlinde algebra: the $T$-plus component  $(\f''(1))^{+}$ 
of $\f''(k=1).$  
Its deformation is $(\f''(1/2-M))^{+}.$  
The corresponding  $q$-deformation
of the {\it big}  Verlinde algebra,   $(\f'(1))^{+},$
is not known.

It is worth mentioning that
there are some integrality/positivity  properties of
the multiplication of the images of the $SL_2$-characters, 
i.e. $p_n^{(1)},$ 
in the Verlinde case.
They will be lost. At least we don't know how to reformulate them for
generic $q.$ Everything else  will survive. See [C5] for
the detail.

Let us calculate the first nontrivial
deformed little Verlinde algebra $V=(\f''(-3/2))^+,$
formed by  $\Q(q^{1/2})$-valued functions on the set 
$\{1/4, 3/4\}.$ The $p$-generators are the images 
$p'_0,\ p'_2$ of the Rogers polynomials 
$$
p_0=1,\ p_2(x) =q^{2x}+q^{-2x}+1+(q+q^{1/2}+1+q^{-1/2}+q^{-1})
$$
for $s=1, \bar{k}=-1/2-1$. The multiplication is standard. We need to know
only
$$
(p'_2)^2=q^{-3/2}\,(1+q^{1/2})^2\, (1+q)^2\,\bigl(p'_2-
q^{-1}\,(1+q^{1/2})\,(1+q^{3/2})\,p'_0\bigr).
$$
The inner product of functions $f,g\in V$ is 
$$
\lan f,g\ran\ =\ fg^*(3/4)+(1-q^{1/2}-q^{-1/2})\, fg^*(1/4).
$$
The Gaussian is proportional to $g(3/4)=1,\ 
g(1/4)=q^{1+\bar{k}}=q^{-1/2}.$

The main formula (the summation with the Gaussian) reads:
$$
\eqalignno{
&p_m p_n(3/4)\, +\, q^{-1/2}(1-q^{1/2}-q^{-1/2})\,p_m p_n(1/4)\ =\cr
&(1-q^{-1})(1+q^{1/2})^{-1}\, q^{{3(m+n)-m^2-n^2\over 4}}\,
p_m(n/2-3/4)\,p_n(3/4).
&\eqnu
\label\verfive\eqnum*
}
$$

Here $q^{1/2}=-\exp(\pi i\om)$, and $0<\om<2/3$ if we want
the inner product to be positive.
Formula (\ref\verfive) describes the Fourier transform
$F_\bullet$ of $\f''(-3/2)$ upon the restriction to $V.$ 
Since $F_0$ fixes $T,$ the restriction is well-defined.
The $V$ is not a $\widetilde{\HH}\ $-module anymore, but
the elements from $\widetilde{\HH}\ $ 
which commute with $T$ act there and make $V$
irreducible. 

The reduction to the corresponding  Verlinde algebra  $V'$
is as follows:
$$
N=5,\ M=2,\ q^{1/2}=-\exp(\pi i/5).
$$ 
Then $p_2(x)$ becomes $q^{2x}+q^{-2x}+1$
and $(p'_2)^2= p'_2+p'_0.$ Also $1-q^{1/2}-q^{-1/2}=(q+q^{-1})^{-2}$
in the inner product.

%
%
%
%
%
%
\AuthorRefNames [BGG]
\references
\vfil

[AI]
\name{R. Askey}, and \name{M.E.H. Ismail},
{ A generalization of ultraspherical polynomials },
in Studies in Pure Mathematics (ed. P. Erd\"os),
Birkh\"auser (1983), 55--78.

[AW]
\name{R. Askey}, and \name{J. Wilson},
{  Some basic hypergeometric orthogonal polynomials
that generalize Jacobi polynomials },
Memoirs AMS {  319} (1985).

[C1]
\name{I. Cherednik},
{ Difference Macdonald-Mehta conjecture },
IMRN {10} (1997), 449--467.

[C2]
\bibline, 
{ Nonsymmetric Macdonald polynomials },
IMRN {10} (1995), 483--515.

[C3]
\bibline
{ On $q$-analogues of Riemann's zeta},
Submitted to Selecta Math. 

[C4]
\bibline, 
{ Nonsymmetric Macdonald polynomials },
IMRN {10} (1995), 483--515.

[C5]
\bibline, 
{ Double Hecke algebras at roots of unity and Gauss-Selberg sums },
In preparation.

[C6]
\bibline
{ From double Hecke algebra to analysis },
Doc.Math.J.DMV, Extra Volume ICM 1998,II, 527--531.

[D]
\name{C.F. Dunkl}, 
{ Hankel transforms associated to finite reflection groups },
Contemp. Math. {138} (1992), 123--138.

[J]
\name{M.F.E. de Jeu},
{ The Dunkl transform }, Invent. Math. {113} (1993), 147--162.

[He]
\name{G.J. Heckman},
{  An elementary approach to the hypergeometric shift operators of
Opdam}, Invent.Math. {  103} (1991), 341--350.

[H]
\name {S. Helgason},
{ Groups and geometric analysis }, 
Academic Press, New York (1984). 

[K]
\name {V.G. Kac},
{ Infinite dimensional Lie algebras },
Cambridge University Press, Cambridge (1990).

[KL]
\name{D. Kazhdan}, and \name{ G. Lusztig},
{ Tensor structures arising from affine Lie algebras. III,}
J. of AMS { 7} (1994), 335--381.

[Ki]
\name{A. Kirillov, Jr.},
{ Inner product on conformal blocks and Macdonald's
polynomials at roots of unity }, Preprint (1995).

[KS]
\name{E. Koelink}, and \name{ J. Stokman},
{ The big $q$-Jacobi function transform },
Publ. IRMA (Univ. Lois Pasteur), 1999/23 (1999).

[M1]
\bibline, {  A new class of symmetric functions },
Publ.I.R.M.A., Strasbourg, Actes 20-e Seminaire Lotharingen,
(1988), 131--171 .

[M2]
\bibline, { Affine Hecke algebras and orthogonal polynomials },
S\'eminaire Bourbaki { 47}:797 (1995), 01--18.

[O1]
\name{E.M. Opdam}, 
{ Some applications of hypergeometric shift
operators }, Invent.Math.{  98} (1989), 1--18.

[O2]
\bibline, { Harmonic analysis for certain representations of
graded Hecke algebras }, 
Acta Math. { 175} (1995), 75--121.

[O3]
\bibline, { Dunkl operators, Bessel functions and the discriminant of
a finite Coxeter group },
Composito Mathematica {85} (1993), 333--373.

\endreferences
\bye